\documentclass[12pt]{amsart}
\usepackage{graphicx}

\newtheorem{theorem}{Theorem}[section]
\newtheorem{prop}[theorem]{Proposition}
\newtheorem{quest}[theorem]{Question}
\newtheorem{lemma}[theorem]{Lemma}

\newtheorem{cor}[theorem]{Corollary}
\newtheorem{defin}[theorem]{Definition}

\newcommand{\map}{\mbox{$\rightarrow$}}
\newcommand{\aaa}{\mbox{$\alpha$}}

\newcommand{\bbb}{\mbox{$\beta$}}

\newcommand{\Ddd}{\mbox{$\Delta$}}
\newcommand{\eee}{\mbox{$\epsilon$}}

\newcommand{\Ggg}{\mbox{$\Gamma$}}
\newcommand{\ggg}{\mbox{$\gamma$}}

\newcommand{\lll}{\mbox{$\lambda$}}

\newcommand{\rrr}{\mbox{$\rho$}}

\newcommand{\Ss}{\mbox{$\Sigma$}}

\newcommand{\ttt}{\mbox{$\tau$}}
\newcommand{\bdd}{\mbox{$\partial$}}

\newcommand{\inter}{\mbox{${\it int}$}}

\begin{document}
\setcounter{tocdepth}{1}

\thispagestyle{empty}
\centerline{\large{{\bf Table of Contents for the Handbook of Knot
Theory}}}

\

\centerline{William W. Menasco and Morwen B. Thistlethwaite, Editors}

\

\begin{enumerate}

\item Colin Adams, {\it Hyperbolic knots}
\item Joan S. Birman and Tara Brendle {\it Braids and knots: A survey}
\item John Etnyre, {\it Legendrian and transversal knots}
\item Greg Friedman {\it  Knot spinning}
\item Jim Hoste {\it The enumeration and classification of knots and
links}
\item Louis Kauffman  {\it Knot Diagrammatics}
\item Charles Livingston  {\it A survey of classical knot concordance}
\item Lee Rudolph  {\it Knot theory of complex plane curves}
\item Martin Scharlemann  {\it Thin position in the theory of 
classical knots }
\tableofcontents
\item Jeff Weeks  {\it Computation of hyperbolic structures in knot
theory}
\end{enumerate}

\pagebreak

  \pagenumbering{arabic}

\title[] {Thin position in the theory of classical knots}

\author{Martin Scharlemann}
\address{\hskip-\parindent
         Mathematics Department\\
         University of California\\
         Santa Barbara, CA 93106\\
         USA}
\email{mgscharl@math.ucsb.edu}

\date{\today}
\thanks{Research supported in part by an NSF grant.}

 %
 %

\maketitle

For our purposes, ``knot theory'' will have the narrowest
interpretation: the study of isotopy classes of locally flat
embeddings of the circle $S^{1}$ in $S^{3}$.  The distinctions between
smooth, PL, and locally flat topological embeddings are usually
unimportant in these dimensions so, for convenience, we'll use
language that is typically associated with the smooth category.  For
example, a knot $K$ is a smooth submanifold of $S^{3}$ diffeomorphic
to $S^{1}$.  The focus of this article will be on a particular
technique for analyzing knots, called ``thin position''.  Much like
the study of crossing diagrams, it's a technique that exploits very
heavily the fact that the ambient manifold is $S^{3}$ and not another
manifold, not even a possibly alternate homotopy sphere.  The roots of
the technique might be traced back to Alexander's proof \cite{Al} of
the Sch\"onfliess theorem, in which he imagines a horizontal plane
sweeping across a sphere embedded in $3$-space and examines the
intersection set during the sweep-out.  The modern use began with a
stunning application by David Gabai, in his proof that knots in
$S^{3}$ satisfy Property R \cite{Ga}.

This is mostly an expository survey; one bit of new mathematics is an
updated application of thin position to the proof that Heegaard
splittings of $S^{3}$ are standard (cf. Corollary \ref{cor:standard}).

\section{From crossing number to bridge number}

One way, perhaps historically the first way, of thinking about knots
in $S^{3}$ is this: Choose a point $p$ in $S^{3}$.  A knot in $S^{3}$
is generically disjoint from $p$, as is an isotopy between knots. 
Thus knot theory in $S^{3}$ (narrowly defined, as above) is equivalent
to knot theory in $\mathbb{R}^{3} \cong S^{3} - \{ p \}$.  Once we
think of a knot as lying in $\mathbb{R}^{3}$, it's natural to imagine
projecting it to $\mathbb{R}^{2}$; this is what we do when we draw the
knot.  A generic such projection will look like an immersed closed
curve in the plane having only double points.  If we keep track, at
the double points, of which strand in the original knot passes over
the other we have the classical description of a knot via its crossing
diagram.  Crossing diagrams are one of the oldest and most naive ways
of trying to classify knots, but their importance has been
re-emphasized by the modern discovery of new knot invariants,
invariants that are most easily described via these projections. 
Unlike more sophisticated invariants (coming, for example, from the
algebraic topology of the knot complement) it is clear that knot
projections make immediate use of the fact that the ambient space is
the $3$-sphere, and not some other manifold.

The first knot invariant that is suggested by knot projection is the
crossing number of the knot.  For a generic projection of the knot to
the plane, the crossing number of the projection is just
the number of double-points of the projection.  An isotopy of the knot
may reduce this number; the crossing number of the knot is defined as
the minimum number that can be achieved via an isotopy of the knot in
$3$-space.  The crossing number is the most natural invariant for
cataloguing knots via their projections (it is the basis of the
standard knot tables) but otherwise it is not a particularly good
invariant.  For a given knot, while it's easy to find an upper bound
for the crossing number (just use any projection), there isn't a
natural way to find the projection that minimizes the crossing number. 
In particular, the behavior of the invariant under e.g. knot sum is
not well-understood.

There is another invariant, much like the crossing number and only
slightly more difficult to describe, that is in fact remarkably well
behaved under knot sum.  Suppose one starts with a knot projection,
and a point on the knot, and starts moving along the knot, recording
at each crossing whether one is on the upper strand (an overcrossing,
say marked with a $+$) or on the lower strand (an undercrossing, say
marked with a $-$).  Continue in this way around the entire knot and
examine the result.  It is a sequence of signs, e.  g. $+, +,
-, +, -, -$.  The number of signs recorded is twice the number of
crossings of the projection since in a trip around the knot, one
passes through each crossing twice.  For the same reason there will be
as many $+$ as $-$ signs.  Now, instead of considering the number of
crossings of the knot (i.  e. focusing on the number of $+$ and $-$
signs) consider the number of times that the sign recorded
{\em changes}.  If the last sign differs from the first, record that
also as a change, as if one were viewing the pattern of $+$ and $-$ on
a circle.  Thus for the sequence $+, +, -, +, -, -$ the number of
changes is $4$.  Of course the number of changes is necessarily even. 
Half that number (i.e. the number of strings of consecutive $+$'s,
say) is called the {\em bridge number} of the knot projection.  Just
as with the crossing number, the bridge number of the projection may
change as the knot is isotoped; the minimum that can be achieved by an
isotopy of the knot is called the bridge number of the knot.  For
reasons that will be apparent, the bridge number of the unknot, which
has no crossings, is conventionally set to be $1$.  This invariant was
introduced by Schubert \cite{Schub}.

The terminology ``bridge number'' is meant to evoke the following
picture: if a knot projection has bridge number $n$, that means there
are $2n$ sign changes from $+$ to $-$ or vice versa; equivalently,
each $+$ sign lies in one of $n$ strings of consecutive $+$ signs. 
Each string of consecutive $+$ signs corresponds to a strand of the
knot which can be thought of as lying just above the plane.  Between
these strands are strands with only $-$ crossings; these can be
thought of as lying on the plane.  Thus the bridge number $\bbb(K)$ of
the knot $K$ is the minimum number of bridges that one would need to
erect on the plane so that the entire knot $K$ could be put on the
plane and on the bridges, with $K$ crossing over each bridge exactly
once.  Schubert \cite{Schub} proved that this invariant is
essentially additive.  To be precise, $\beta - 1$ is additive; that
is, $$\beta(K_{1} \# K_{2}) - 1 = (\beta(K_{1}) - 1) + (\beta(K_{2}) -
1)$$ or $$\beta(K_{1} \# K_{2}) = \beta(K_{1}) + \beta(K_{2}) - 1.$$ A
modern proof is given in \cite{Schul}.

Schubert's remarkable result suggests that there should be a more
natural way of viewing bridge number.  Consider the description of the
knot just given: lying mostly on the plane, but with certain sections
of it elevated above the plane on bridges.  The original perspective
on the knot is the bird's-eye view in which we look down on the knot
from above.  Imagine instead the perspective of someone standing on
the plane, looking sideways at the knot and thinking about the height
of the knot above the plane.  That is, instead of the projection
$\mathbb{R}^{3} \map \mathbb{R}^{2}$ that described the original knot projection,
consider instead the projection $h:\mathbb{R}^{3} \map \mathbb{R}$
defined by height above the plane.  Then $h|K$ will achieve a single
maximum on each bridge, and each pair of successive maxima is
separated by a minimum corresponding to a strand of $K$ that lies on
the plane.  In other words, if $K$ is put on $n$
bridges, then $h|K$ has $n$ maxima.  This number of maxima can be
preserved even when $K$ is moved slightly to make it generic with
respect to the height function by, for example, putting a bit of a dip
into the level strands lying on the plane.  In other words, we have
the following:

\begin{prop} Suppose $K \subset \mathbb{R}^{3}$ has a projection with
bridge number $n$.  Then $K$ may be isotoped so that the standard
height function $h: \mathbb{R}^{3} \map \mathbb{R}$ restricts to a
Morse function on $K$ with $n$ maxima and $n$ minima.
\end{prop}

Conversely, we have

\begin{prop} Suppose $K \subset \mathbb{R}^{3}$, and the standard
height function $h: \mathbb{R}^{3} \map \mathbb{R}$ restricts to a
Morse function on $K$ with $n$ maxima and $n$ minima.  Then $\beta(K)
\leq n$.
\end{prop}

\begin{figure} [tbh]
\centering
\includegraphics[width=0.5\textwidth]{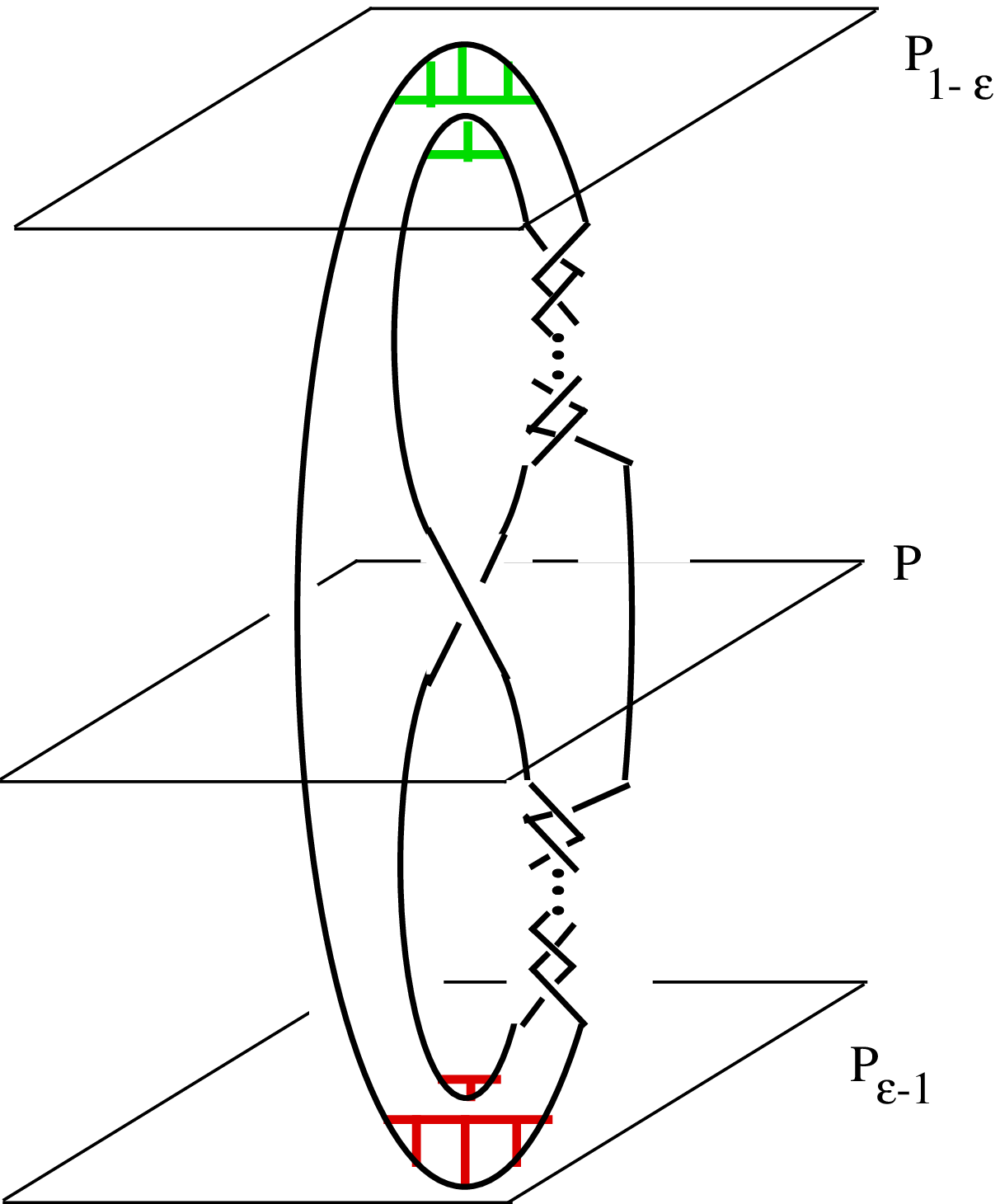}
\caption{} \label{fig:bridge}
\end{figure}

\begin{proof} Each point on $K$ at which $h|K$ has a local maximum can
be pushed even higher by an isotopy of $K$, say along a path that
rises from the maximum and, by general position, misses the rest of
$K$ as it rises.  Similarly, each point at which $h|K$ has a local
minimum can be pushed lower.  So with no loss, $K$ may be isotoped so
that all the maxima occur near the same height (say $1$) and similarly
all the minima occur near the same height (say $-1$).  Consider how the 
plane $P = h^{-1}(0)$ divides $K$:  $P$ cuts $K$ into $2n$ strands, 
$n$ above the plane, each containing a single maximum, and $n$ below 
the plane, each having a single minimum.  

The plane $P_{\eee - 1}$ slightly above height $-1$ cuts off $n$ tiny
strands of $K$, one for each minimum, that lie below the plane.  These
are clearly isotopic rel their endpoints to $n$ disjoint arcs in
$P_{\eee - 1}$.  See Figure \ref{fig:bridge}.  Between the heights
$\eee - 1$ and $0$ there are no critical points of $h$ on $K$.  It
follows that in fact all $n$ strands of $K$ lying below $P$ can be
simultaneously isotoped rel their ends to lie on $P$.  Similarly, all
strands of $K$ lying above $P$ can be isotoped rel their end points so
that they consist of level arcs in a plane $P_{\eee}$ just above $P$,
together with vertical arcs at their end points between $P$ and
$P_{\eee}$.  When viewed from above as projected onto $P$, $K$ then
has a projection with bridge number no higher than $n$.  \end{proof}

Combining the two propositions above gives a more natural definition
of bridge number of a knot $K$: take the number of maxima that the
standard height function $h$ has on the knot (in general position with
respect to $h$) and minimize that number via an isotopy of $K$.  The
result is the bridge number $\beta(K)$ of $K$.  If there is a
horizontal plane $P$ with the property that all maxima of $h|K$ lie
above $P$ and all minima lie below it, then $K$ is said to be in {\em
bridge position} with respect to $h$; the plane $P$ is called a {\em
dividing plane} for $K$.  Any knot $K$ can be isotoped so that it is
in bridge position, with $\beta(K)$ maxima and $\beta(K)$ minima.

\section{From bridge number to width}

Consider the critical values of $K \subset \mathbb{R}^{3}$ in general
position with respect to $h: \mathbb{R}^{3} \map \mathbb{R}$, the
standard height function.  As noted above, it's always possible to
isotope $K$ so as to raise the height of a maximum or to lower the
height of a minimum, without affecting the height of any other
critical point.  Indeed there is such an isotopy whose support on $K$
is limited to a neighborhood of the critical point whose height is
changed.  Similarly, there is no difficulty lowering the height of a
maximum through an interval that contains no other critical values,
since by standard Morse theory, the preimage of an interval without
critical values is a simple product.  Combining these two
observations, it's easy to isotope $K$ to interchange the heights of
critical points whose critical values are adjacent, so long as both
critical points are maxima, or both are minima.

\begin{figure} [tbh]
\centering
\includegraphics[width=0.2\textwidth]{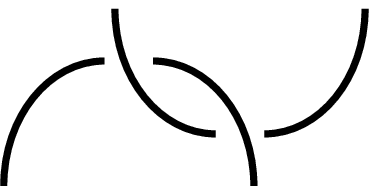}
\caption{} \label{fig:stuck}
\end{figure}

The only difficulty in rearranging the heights of critical points is
the interchange of two adjacent critical values in which the higher
value is a maximum and the lower value is a minimum.  Informally, this
can be described as moving a maximum down past a minimum.  Such a move
may or may not be possible, depending on the structure of the knot (cf
Figure \ref{fig:stuck}).  It's reasonable then to think of this move
(pushing a maximum down past a minimum) as simplifying the picture of
the knot, when it can be done; the point of thin position is to
formally capture this idea in a useful way.

\begin{figure} [tbh]
\centering
\includegraphics[width=0.9\textwidth]{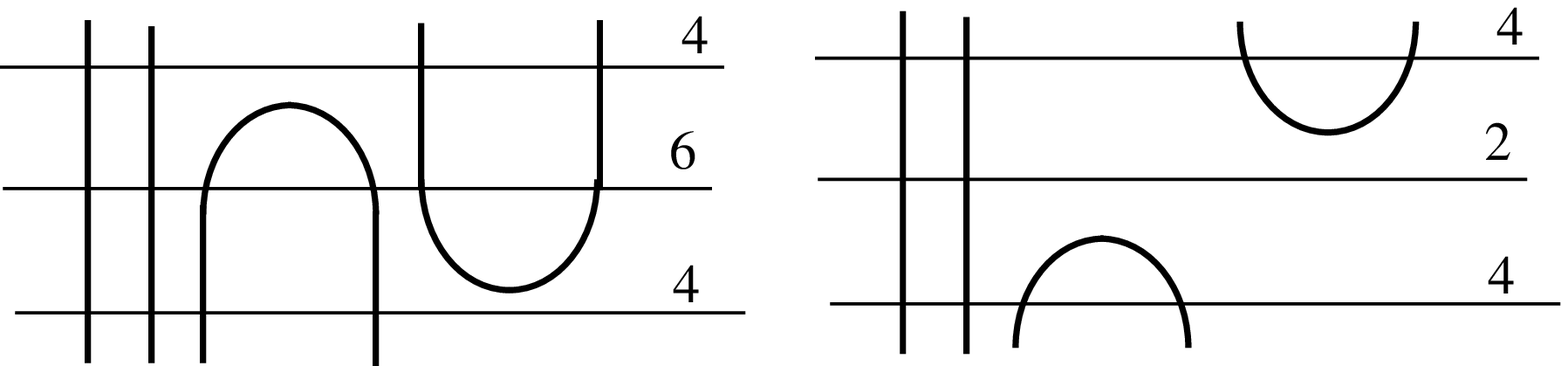}
\caption{} \label{fig:thin}
\end{figure}

Let $K \subset \mathbb{R}^{3}$ be a knot in general position with
respect to the standard height function $h: \mathbb{R}^{3} \map
\mathbb{R}$.  That is, $h|K$ is a Morse function for which no two
critical points have the same critical value.  For each $t \in R$ let
$P_{t}$ denote the plane $h^{-1}(t)$.  If $t$ is a regular value for
$h|K$ then $K$ crosses $P_{t}$ transversally, necessarily in an even
number $w(t) \in \mathbb{N}$ of points.  The number $w(t)$ changes
only at critical values, where it increases by two at each minimum and
decreases by two at each maximum.  So if $c_0 < c_{1} < \dots < c_n$
are the critical values of $h|K$ and values $r_1, \dots, r_{n}$ are
chosen so that $c_{i-1} < r_i < c_{i},$ $ i = 1, \ldots, n$, then the
function $w(t)$ is determined by the sequence $w(r_{i}),$ $ i = 1,
\ldots n$.  That sequence is unchanged by pushing a maximum down past
another maximum (or a minimum past another minimum) but is affected by
pushing a maximum down past a minimum.  In the last case, if the
critical values are $c_{i-1}, c_{i}$, then the reordering changes
$w(r_{i})$ to $w(r_{i}) - 4$ and has no other effect.  See Figure
\ref{fig:thin}. More dramatically, if $K$ can be isotoped so that the
maximum and its adjacent minimum cancel, then both $c_{i-1}$ and
$c_{i}$ disappear and so both $w(r_{i}), w(r_{i+1})$ disappear from
the sequence.

These
changes are so straightforward, they suggest the following definition:

\begin{defin} \label{defin:width}
{\rm Suppose, as above, $K \subset \mathbb{R}^{3}$ is in general position with
respect to the standard height function $h$, $c_0 < c_{1} < \dots <
c_n$ are the critical values of $h|K$ and regular values $r_{i} \in R$
are chosen so that $c_{i-1} < r_i < c_{i}, i = 1, \ldots, n$.  The
{\em width of K with respect to $h$}, denoted by $w(K, h)$, is $\sum_i
w(r_{i})$.  The {\em width of K}, denoted by $w(K)$, is the minimum of
$w(K', h)$ over all knots $K'$ isotopic to $K$.  We say that $K$ is in
{\em thin position} if $w(K, h) = w(K).$}
\end{defin}

There is a small technical advantage in modifying this definition
slightly so that it can be applied directly and a bit more usefully in
the compact situation $K \subset S^{3}$.  Define the standard height
function $h: S^{3} \map \mathbb{R}$ as the composition $S^{3} \subset
\mathbb{R}^{4} \map \mathbb{R}$ of the inclusion and the standard
projection onto the last factor of $\mathbb{R}^{4}$.  Then $h$
has two critical points on $S^{3}$ (called the north and south poles). 
A knot $K$ in general position with respect to $h$ will be disjoint
from these poles, and the above definition can be restated for this
height function.  The most significant difference with this definition
is that (for $-1 < t < 1$), $h^{-1}(t)$ is now a $2$-sphere $P_{t}$
instead of a plane; if $K \subset S^{3}$ is in bridge position, then
what divides the maxima from the minima is a {\em dividing sphere}
instead of a dividing plane.  Whenever $t$ is a regular value for
$h|K$, we continue to denote $|P_{t} \cap K|$ by $w(t)$.

Note that if $K \subset \mathbb{R}^{3}$ or $K \subset S^{3}$ has been isotoped
to be in thin position with respect to the standard height function
$h$, it is impossible to have a limited isotopy that simply pushes a
maximum below a minimum or cancels a maximum with a minimum, since
either move decreases the width by at least $4$.

\section{Application: Thinning the unknot}

To show the power of this idea, we begin with a simple exercise that
illustrates how thin position interacts with geometric properties of a
knot, in particular with an essential surface in the knot complement.

It will be useful to have the following notation and definitions: 

\bigskip

\noindent {\bf Notation:} For $M$ a manifold and $X \subset M$ a
polyhedron, let $\eta(X)$ denote a closed regular neighborhood,
whereas (abusing notation slightly) $M - \eta(X)$ will mean the closed
complement of $\eta(X)$ in $M$.

\bigskip

\begin{defin} 
{\rm Let $K \subset S^{3}$ be a knot, $P \subset S^{3}$ be a sphere that is
level with respect to the standard height function and is transverse
to $K$.  Let $B_u$ and $B_l$ denote the balls which are the closures
of the region above $P$ and below $P$ respectively.  An {\em upper
disk} (resp.  {\em lower disk}) with respect to $P$ is a disk $D
\subset S^3 - \eta(K)$ transverse to $P$ such that $\bdd D = \alpha
\cup \beta$, where $\beta$ is an arc imbedded on $\partial\eta(K)$,
parallel to a subarc of $K$, $\alpha$ is an arc properly imbedded in
$P - \eta(K)$, $\bdd \alpha = \bdd \beta$ and a small product
neighborhood of $\aaa$ in $D$ lies in $B_u$ (resp.  $B_l$) i.e., it
lies {\em above} (resp.  {\em below}) $P$.}
\end{defin}

Note that $\inter(D)$ may intersect $P$ in simple closed curves or
indeed in other arcs.  An innermost simple closed curve cuts off a
disk that lies either entirely above or below $P$.  Such a disk is
called respectively an {\em upper cap} or {\em lower cap}.

Natural upper disks with interiors disjoint from $P$ arise, for
example, in the case where the arcs $K \cap B_{u}$ each have exactly
one maximum, so in particular the collection of arcs is the untangle. 
To see the upper disks, consider what happens as a descending level
sphere $P_{t}$ sweeps across a maximum.  Join the descending arcs of
$K$ from that maximum by an arc in $P_{t}$ and continue to carry that
arc down all the way to $P$.  The result is a disjoint collection of
upper disks, one for each component of $K \cap B_{u}$; each
intersects a level sphere in at most one arc.  These are called a
family of {\em descending disks} for $K \cap B_{u}$.  Choices are
involved in this construction: each time a new maximum is encountered
as $P_{t}$ sweeps down, a choice is made about how existing arcs from
earlier descending disks lie in relation to this maximum.  In fact, if
one allows isotopies that raise and lower maxima (but never introduce
minima) one has the general observation, whose proof is mostly left as
an exercise:

\begin{lemma} \label{lemma:descending} Suppose $P$ is a level sphere
for $K \subset S^{3}$ and each component of $K \cap B_{u}$ has a
single maximum.  Suppose $\Ddd$ is a collection of disjoint upper
disks contained entirely in $B_{u}$.  Then $\Ddd$ can be isotoped rel
the arcs $\Ddd \cap P$ so that it becomes part of a complete
collection of descending disks.  Moreover, such an isotopy can be
found so that during the isotopy no new critical points of $h|K$ are
introduced.
\end{lemma}

A proof hint is this: Start with any complete collection of descending
disks $\Ddd'$ and alter $\Ddd$ and $\Ddd'$ to reduce $|\Ddd' \cap
\Ddd|$, i.e. the number of components in $\Ddd' \cap \Ddd$.

\begin{prop} \label{prop:thinbridge} Suppose $K \subset S^{3}$ is the
unknot, in bridge position with respect to the standard height
function $h: S^{3} \map \mathbb{R}$.  There is a dividing sphere $P$
for $K$ so that a maximum and a minimum can simultaneously be
isotoped to lie on $P$.  During the generic isotopy, the width
remains unchanged.
\end{prop}

\begin{proof}  Since $K$ is the unknot, it bounds a disk $D$.  By a 
small isotopy of $D$ we can arrange that $D$ is in general position 
with respect to $h$ and near each critical point of $h|K$, $h|D$ has 
a half-center singularity. (See Figure \ref{fig:saddlecenter}, applied 
there in another context.) That is, near each maximum of $K$, $D$ is 
incident to $K$ from below and near each minimim, $D$ is incident to 
$K$ from above.  In particular, for a level sphere $P$ just below a 
maximum (resp. just above a minimum) one of the components of $D - P$ 
is an upper disk (resp. lower disk) contained entirely above (resp. 
below) $P$.  

Let $t_{l} < t_{u}$ be the heights of, respectively, the highest
minimum and the lowest maximum of $h|K$ and for each $t_{l} < t <
t_{u}$, let $P_{t}$, as above, denote the level sphere $h^{-1}(t)$. 
We have just seen that for $t$ slightly less than $t_{u}$, $D - P_{t}$
contains an upper disk among its components and similarly, for $t$
slightly greater than $t_{l}$, $D - P_{t}$ contains a lower disk among
its components.

{\bf Claim:} There is a value of $t, t_{l} < t < t_{u}$ for which
$P_{t}$ admits disjoint upper and lower disks $D_{u}, D_{l}$ so that no
component of $\inter(D_{u}) \cap P$ or $\inter(D_{l}) \cap P$ is an arc.

{\bf Proof of claim:} For any generic value of $t, t_{l} < t < t_{u}$
(i.  e. a value of $t$ for which $P_{t}$ is transverse to $D$),
consider an outermost arc of $P_{t} \cap D$ in $D$.  The disk it cuts
off from $D$ is either an upper disk or a lower disk; moreover the
subarc of $K$ incident to the disk lies either entirely above or
entirely below $P_{t}$.  Hence for each generic value there is either
an upper or a lower disk as desired.

Now imagine $t$ ascending from $t_{l}$ up to $t_{u}$.  Since near
$t_{l}$ an outermost arc of $P_{t}$ in $D$ cuts off a lower disk and
near $t_{u}$ one cuts off an upper disk, and at any generic $t$ in
between, one or the other is cut off, there are two possibilities. 
One is that there is a generic value for which outermost arcs of $P
\cap D$ in $D$ cut off both an upper and a (disjoint) lower disk from
$D$; then we are done with proving the claim.  The second possibility
is that there is a critical value $t_{0}, t_{l} < t_{0} < t_{u}$ for
$h|D$, whose critical point is necessarily an interior point of $D$,
so that for small $\eee$, outermost arcs of $P_{t_{0} + \eee}$ and
$P_{t_{0} - \eee}$ in $D$ cut off respectively an upper and a lower
disk.  In this case, thicken $D$ slightly and let $D_{\pm}$ be the
boundary disks of the thickened region.  That is, $D_{\pm}$ are two
copies of $D$, very near to $D$ but on opposite sides.  Then
$P_{t_{0}}$ is transverse to both $D_{+}$ and $D_{-}$ and outermost
arcs of $P_{t_{0}}$ in $D_{+}$ and $D_{-}$ will cut off
(disjoint) upper and lower disks, completing the proof of the claim.

\bigskip

In the special case in which both upper and lower disks have interiors 
entirely disjoint from $P$, Lemma \ref{lemma:descending} applies, and 
$K$ may be isotoped rel $K \cap P$, never changing the width, so that 
afterwards the upper disk is a descending disk and, dually, the lower 
disk is an ascending disk.  These disks then define the isotopy of a 
maximum and minimum to disjoint arcs in $P$, as required to complete 
the proof of the Proposition.

If the interiors of the upper disk $D_{u}$ or the lower disk $D_{l}$ 
are not disjoint from $P$, the argument is only moderately more 
complicated.  In that case, each component of intersection is a closed 
curve, and an innermost such closed curve on $D_{u}$ or $D_{l}$ cuts 
off a disk that is an upper or lower cap.  Suppose, for example, that 
there is an upper cap $C$.  A standard innermost disk, outermost 
arc argument will alter a complete collection of descending disks for 
$K \cap B_{u}$ to a collection of upper disks disjoint from $C$.  Via 
Lemma \ref{lemma:descending} there is an isotopy of $K \cap B_{u}$ rel 
$P$ which does not increase width so that afterwards, these disks are 
descending disks, possibly now again intersecting $C$, but only in 
their interiors.  Then alter $C$, via an innermost disk argument, 
isotoping $C$ so that afterwards the set of descending disks is 
disjoint from $C$.  This establishes that, after an isotopy of $K$ and 
$C$ with support away from $P$ and never increasing the width, there 
is a complete collection of descending disks for $K \cap B_{u}$ that 
is disjoint from $C$.

Suppose then that there are disjoint upper and lower caps $C_{u}$ and
$C_{l}$.  Modify $K \cap B_{u}, K \cap B_{l}$ away from $P$ so that
all descending disks (resp.  all ascending disks) are disjoint from
$C_{u}$ (resp.  $C_{l}$).  $\bdd C_{u}$ and $\bdd C_{l}$ bound
disjoint disks $E_{u}, E_{l}$ in the sphere $P$.  Pick a component of
$K \cap B_{u}$ incident to $E_{u}$ and a component of $K \cap B_{l}$
incident to $E_{l}$.  A descending disk for the former will intersect
$P$ inside of $P_{u}$ and an ascending disk for the latter will
intersect $P$ inside of $P_{l}$.  In particular, they will be 
disjoint, and so they can be used to isotope arcs to $P$ as required.
A similar argument applies if there is an upper cap and a lower disk 
whose interior is disjoint from $P$, or symmetrically.  

The only remaining case to consider is when there are no upper caps,
say, but the interior of the upper disk $D_{u}$ intersects $P$, so
there are lower caps.  We will show that in this case, one of the 
cases we have already considered also applies.

Let $\Ddd$ be a complete collection of descending disks for $K \cap
B_{u}$.  We argue by induction on $|D_u \cap \Ddd|$ that there are
both an upper disk (or an upper cap) and a lower cap so that their
boundaries are disjoint.  If $|D_u \cap \Ddd| = 0$ then each component
of $D_u \cap B_u$ lies in the ball $B_{u} - \Ddd$.  After compressing
in $B_{u} - \Ddd$ each component becomes a disk.  Since a neighborhood
of $\bdd D_u$ lies in $B_u$, it follows that after the compressions,
$\bdd D_{u}$ bounds a disk in $B_u - \Ddd$, hence in $B_{u} - K$, as
required.

So suppose $D_u \cap \Ddd \neq \emptyset$.  A simple innermost disk
argument could eliminate a closed curve of intersection, so we can
take all components of intersection to be arcs.  Surprisingly, we may
also assume that the lower cap is a slight push-off of a disk
component of $D_u \cap B_l$.  Indeed, consider an innermost disk of
$D_u - P$.  If it lies in $B_u$ then it is an upper cap disjoint from
the lower cap and we are done.  If it lies in $B_l$ then we may as
well take a slight push-off as our lower cap.

This surprising fact means that an outermost arc of $D_u \cap
\Ddd$ in $\Ddd$ can be used to $\bdd$-compress $D_u \cap B_u$ to an
arc that is disjoint from the lower cap.  This boundary compression
defines an isotopy on the interior of $D_u$ that reduces $|D_u
\cap \Ddd|$ without disturbing the disjoint lower cap.  
After the isotopy, the result follows by induction.
\end{proof}

Note that if the ends of the maximum and minimum arcs given by Proposition 
\ref{prop:thinbridge} both coincide then they constitute all of $K$ 
and $K$ already was in thin position.  Otherwise, the minimum and 
maximum can be pushed on past each other, or just cancelled if the arcs 
have a single end in common, reducing the width.  Thus we have:

\begin{cor} \label{cor:thinbridge} If the unknot is in bridge
position, then either it is in thin position (and so has just a single
minimum and maximum) or it may be made thinner via an isotopy that
does not raise the width.
\end{cor}

The corollary begs the question: is the hypothesis that the unknot is
in bridge position really needed?  That is

\begin{quest} \label{quest:unknot} Suppose $K \subset S^{3}$ is the
unknot.  Is there an isotopy of $K$ to thin position (i.  e. a single
minimum and maximum) via an isotopy during which the width is
never increasing?
\end{quest}

One suspects there are counterexamples, though it would be difficult 
to prove for such a counterexample that no such isotopy exists.

\section{Thick and thin regions}

If $K \subset S^{3}$ is in bridge position, then $w(t), -1 < t < 1$ is
constant on intervals that contain no critical values of $h|K$. 
$w(t)$ always increases by $2$ at each critical value $h|K$ that lies
below the height of a dividing sphere and then decreases by $2$ at
each critical value that lies above the height of a dividing sphere. 
If $K$ is not in bridge position, $w(t)$ will still increase or
decrease by $2$ at each critical value of $h|K$, but $w(t)$ will have
one or more local minima as well as more than one local maximum.  For
example, if $t_{0}$ is a regular value of $h|K$ and the critical
values above and below $t_{0}$ correspond respectively to a minimum and
maximum of $h|K$, then $w(t)$ will be greater if $t$ is either
increased or decreased past the adjacent critical values of $h|K$. 
That is, $w(t_{0})$ is a local minimum of $w(t)$; the interval of
regular values for $h$ on which it lies is called a thin region (and
the corresponding heights the thin levels).  The level sphere
$P_{t_{0}}$ is called a thin sphere.

Symmetrically, if $t_{0}$ is a regular value of $h|K$ and the critical
values above and below correspond respectively to maxima and minima of
$h|K$, then the interval of regular values for $h$ on which it lies is
called a thick region (and the corresponding heights the thick
levels).  The level sphere $P_{t_{0}}$ is then called a thick sphere.

To repeat in the notation of Definition \ref{defin:width}, call the
level $r_i$ a {\em thin level} of $K$ with respect to $h$ if the
critical point $c_i$ is a local maximum for $h$ and $c_{i+1}$ is a
local minimum for $h$.  Dually $r_i$ is a {\em thick level} of $K$
with respect to $h$ if $c_i$ is a local minimum and $c_{i+1}$ is a
local maximum.  Many $r_{i}$ may be neither thin nor thick.  Since the
lowest critical point of $h|K$ is a minimum and the highest is a
maximum, there is one more thick level than thin level.

There is an alternative way, using thin and thick levels, to calculate
the width $w(K, h)$ of a knot in $S^{3}$.  Choose values $r_{i_{1}}
\dots, r_{i_k}$ in Definition \ref{defin:width} to be those of the
thick levels of $K$ and $r_{j_1}, \dots, r_{j_{k-1}}$ to be those of
the thin levels, so $r_{i_{l}} < r_{j_{l}} < r_{i_{l+1}}, 1 \leq l \leq
k-1$.

\begin{lemma}
Let $a_{l} = w(r_{i_{l}})$ and $b_{l} = w(r_{j_{l}})$.  Then \[2w(K) = 
\sum_{l = 1}^k a_{l}^2 - \sum_{l = 1}^{k-1} b_{l}^2.\]
\end{lemma} 

\begin{proof}
See the last section of \cite{ScSc} for McCrory's simple proof; indeed
a contemplative look at the last figure there should suffice. 
\end{proof}

If $K$ is in thin position, thin and thick levels of the height
function have important geometric properties.  For example, Thompson
\cite{Th1} showed that if thin position for $K$ is not bridge
position, so $K$ has thin levels, then there is an essential meridinal
planar surface for $K$.  One way of finding such a surface was
recently identified by Ying-Qing Wu \cite{Wu}:

\begin{theorem} \label{theorem:wu} Suppose $K \subset S^{3}$ is in
thin position but not in bridge position, so there are thin levels. 
Suppose $P_{r_{j_{l}}}$ is the thinnest thin sphere (that is, among
all values at thin levels $r_{j_{i}}$, $w(r_{j_{l}})$ is the lowest). 
Then the planar surface $P_{-} = P_{r_{j_{l}}} - \eta(K)$ is essential
in $S^{3} - \eta(K)$.  That is, $P_{-}$ is incompressible in $S^{3} -
\eta(K)$ and is not a boundary parallel annulus.
\end{theorem} 

A sample application for this result comes from work of Gordon and
Reid \cite{GR}.  A knot $K \subset S^{3}$ is said to have {\em tunnel
number one} if there is a properly imbedded arc $\ttt \subset S^{3} -
\eta(K)$ so that $S^{3} - (\eta(K) \cup \eta(\ttt))$ is a genus two
handlebody.  Gordon and Reid showed that a tunnel number one knot has
no incompressible planar surfaces in its complement.  Combining the 
two, we have Thompson's result \cite{Th1}:

\begin{theorem} Suppose a tunnel number one knot $K \subset
S^{3}$ is in thin position.  Then it is also in minimal bridge position.
\end{theorem}

A second feature of thin and thick spheres is that they intersect
essential surfaces in the knot complement in a controlled way.  For
example, suppose $F$ is a Seifert surface for $K$, i.  e. an
orientable surface in $S^{3}$ whose boundary is $K$.  Suppose $F$ is
in general position with respect to the height function $h: S^{3} \map
\mathbb{R}$.  That is, all the critical points of $h|F$ and $h|\bdd F = K$ are
non-degenerate and no two occur at the same height.  First consider
thin levels:

\begin{theorem} \label{theorem:arcthin} Suppose $K \subset S^{3}$ is
in thin position but not in bridge position, so there are thin levels. 
Suppose $F$ is a Seifert surface for $K$, in general position with
respect to $h$, and $P_{r}$ is a thin sphere.  Then every arc
component of $F \cap P_{r}$ is essential in $F$.
\end{theorem}

\begin{proof} The argument is most easily described in the compact
manifold $S^{3} - \eta(K)$ so let $P^{-}$ be the planar surface $P_{r}
- \eta(K)$.  Since $F$ is a Seifert surface we can assume that $F \cap
\bdd\eta(K) = \bdd F$ is a longitude, and so $\bdd F$ intersects each
component of $\bdd P^{-}$ exactly once.  In particular, every arc
component of $F \cap P^{-}$ is essential in $P^{-}$; indeed such
components pair up the components of $\bdd P^{-}$.  Suppose some
arc component is inessential (i.  e. $\bdd$-parallel) in $F$; let $\aaa$
be an outermost such component, i.  e. a component cutting off a
subdisk $E$ of $F$ which contains no other arc component of $F \cap
P^{-}$, though it may contain circle components.  

\begin{quotation}

{\bf Aside:} An experienced $3$-manifold topologist would expect first
to eliminate these circle components, but in fact we do not know that
we can, for although Theorem \ref{theorem:wu} tells us that the
thinnest level sphere gives rise to an incompressible surface, we do
not know this to be true for an arbitrary thin level sphere and so we
cannot automatically eliminate a circle in $F \cap P_{-}$ just because
it bounds a disk in $F - P_{-}$.

\end{quotation}

Now $\bdd E$ consists of two arc components: $\aaa \subset P_{-}$ and
a subarc $\bbb$ of $\bdd F = K$.  Hence the disk $E$ can be used to
isotope $\bbb$ to $\aaa$, though note that this isotopy may move
$\bbb$ through $P_{-}$ since the interior of $E$ may have circles of
intersection with $P_{-}$.  Nonetheless, we do know that $\bbb$ lies
entirely on one side of $P_{-}$, say above $P_{-}$.  So the effect of
moving $\bbb$ to $\aaa$ is at least the effect of moving a maximum
(namely a maximum of $\bbb$) past a minimum (namely the minimum that
is the lowest critical value of $h|K$ above height $r$).  In fact, 
much more may be accomplished, e. g. the elimination of other critical 
points from $h|\bbb$, but the net effect is to lower the width of 
$K$.  Since we have assumed that $K$ begins in thin position, this is 
impossible, proving the theorem.
\end{proof}

Note that the fact that the isotopy of $\bbb$ may pass through $P_{-}$ 
is an example of why this argument cannot be directly applied to 
Question \ref{quest:unknot}: we can always thin the unknot via an arc 
by arc isotopy as just described, but we have little control over the 
width during each of these isotopies. 

It is a bit more surprising that there is a version of Theorem
\ref{theorem:arcthin} that also applies at a thick level:

\begin{theorem} \label{theorem:arcthick} Suppose $K \subset S^{3}$ is
in thin position and $F$ is a Seifert surface for $K$ in general
position with respect to $h$.  Suppose $P_{r}$ is a thick sphere for
$K$.  Let $c_{-} < r < c_{+}$ be critical values of $h|K$ that are
adjacent to $r$, so in particular $c_{+}$ is the height of a maximum
and $c_{-}$ is the height of a minimum.  Then either $K$ is the unknot
or there is a value $r', c_{-} < r' < c_{+}$ so that every arc
component of $F \cap P_{r'}$ is essential in $F$.
\end{theorem}

\begin{proof} Consider a level sphere $P_{+}$ just below height
$c_{+}$, so in particular there are no critical values for $F$ between
the level of $P_{+}$ and $c_{+}$.  If all arc components of $F \cap
P_{+}$ are essential in $F$, we are done, so suppose $E \subset F$ is
a disk cut off by an outermost inessential arc $\aaa$, with $\bdd E =
\aaa \cup \bbb,\bbb \subset K$.  Let $\ggg \subset (K - P_{+})$ be the
interval containing the maximum at height $c_{+}$.  Suppose first that
$\bbb$ lies below $P_{+}$, so in particular $E$ is a lower disk.  If
the ends of $\bbb$ coincide with the ends of $\ggg$ then $E$ 
describes an isotopy of $\bbb$ up to $P_{+}$; after the isotopy $K$ 
has a single maximum and minimum and so is the unknot.  If a single 
end of $\bbb$ coincides with a single end of $\ggg$ then $E$ can be 
used to isotope $\bbb$ up past $c_{+}$, cancelling the maximum in 
$\ggg$, as well as one or more critical points in $\bbb$.  This would 
reduce the width of $K$, which is impossible.  Finally, if the ends 
of $\bbb$ and $\ggg$ are disjoint, then $E$ can be used to move 
$\bbb$ above $c_{+}$ (since there are no critical values of $h|F$ 
between the level of $P_{+}$ and $c_{+}$) thereby moving a minimum 
past a maximum, and possibly cancelling other critical points on 
$\bbb$.  Again this would contradict the assumption that $K$ is in 
thin position.   We therefore conclude that in fact $\bbb$ lies above 
$P_{+}$ (indeed perhaps $\bbb = \ggg$) so $E$ is an upper disk.

Similarly, a level sphere $P_{-}$ just above the level of $c_{-}$
either cuts off a lower disk or we are done.  If there is a generic
height between $P_{\pm}$ for which all arcs of intersection with $F$
with the corresponding level sphere are essential in $F$ we are done. 
On the other hand, if at every generic level there is at least one
inessential arc of intersection then there is always a disk cut off
from $F$ that is an upper or a lower disk.  Then, as in the proof of
Proposition \ref{prop:thinbridge} (perhaps after thickening $F$ as
there we thickened $D$), there is a level sphere $P$ that cuts off
simultaneously an upper disk $E_{u}$ and a disjoint lower disk
$E_{l}$, via arcs $\aaa_{u}, \aaa_{l} \subset (P \cap F)$.  If both
ends of $\aaa_{u}$ and $\aaa_{l}$ coincide, then $K$ is the unknot. 
Otherwise, $E_{u}$ and $E_{l}$ can be used to isotope arcs of $K$ to
$P$, lowering the width as described in the proof of Proposition
\ref{prop:thinbridge}.  Since $K$ was assumed to be thin, this is 
impossible.    \end{proof}

The original application by Gabai that prompted the definition of thin
position is in a similar setting \cite{Ga}.  Gabai's application was
in the proof of the Poenaru conjecture: For $K \subset S^{3}$ there is
an essential (e.  g. non-separating) planar surface $(Q, \bdd Q)
\subset (S^{3} - \eta(K), \bdd \eta(K))$ whose boundary components are
all longitudes (i.e. of slope $0$) on $\bdd \eta(K)$ only if $K$ is
the unknot.  For technical reasons, Gabai wished to exhibit a planar
surface $(P, \bdd P) \subset (S^{3} - \eta(K), \bdd \eta(K))$
transverse to $Q$ for which $\bdd P$ is a collection of meridians of
$\eta(K)$ and every arc of $P \cap Q$ is essential in both $P$ and
$Q$.  Gabai applied essentially the argument above, substituting $Q$
for $F$.  The upshot is a pair of planar surfaces, $P, Q \subset S^{3}
- \eta(K)$ so that in each surface, viewing the boundary components as
large vertices and the intersection arcs as edges connecting the
vertices, we have what appears to be a planar graph, with all the rich
structure that this implies.  (It must be said, though, that this
structure is only the starting point of Gabai's deep and complex
argument using sutured manifold theory.)  This application is an echo,
in some sense, of Laudenbach's \cite{Lau} seminal introduction of graph
theory into such arguments.  Laudenbach proved the Poenaru conjecture,
even for knots in a mere homotopy $3$-sphere, in the simplest
interesting case: when $Q$ has just three boundary components.

\bigskip

There are two important directions in which Gabai's application
generalizes.  To formulate the first extension, note first that the
Poenaru conjecture implies in particular that if $0$-framed surgery on
a knot $K$ yields a manifold containing a non-separating sphere, then
$K$ itself has genus $0$.  Gabai generalizes this to show that if
$0$-framed surgery on $K$ yields a manifold containing a
non-separating genus $g$ surface, then the genus of $K$ is no larger
than $g$.  The relevant extension of the thin-position part of his
argument is to replace the punctured sphere $Q$ with an essential
punctured genus $g$ surface; $P$ remains a planar surface intersecting
$Q$ in arcs that are essential in both surfaces. 

A second extension is used in the celebrated proof by Gordon and
Luecke \cite{GL}, that a knot is determined by its complement.  In
their setting, surgery on $K$ with slope $\pm 1$ hypothetically yields
$S^{3}$ again, with the core of the solid torus representing a new
knot $K' \subset S^{3}$.  They seek to simultaneously find meridinal
planar surface $P$ and $Q$ in $S^{3} - \eta(K)$ and $S^{3} - \eta(K')$
respectively, so that $P$ and $Q$ are transverse and all arcs in $P
\cap Q$ are essential in both surfaces.  Viewed in $S^{3} - \eta(K)$,
each component of $\bdd Q$ has slope the same as the surgery slope. 
In order to find the pair of surfaces, they must simultaneously use
thin position on height functions for $K \subset S^{3}$ and for $K'
\subset S^{3}$.

It's important to repeat that in both the Gabai argument and the
Gordon-Luecke argument, the use of thin position is only one of many
parts of the full proof, and far from the deepest.

\bigskip

\begin{figure} [tbh]
\centering
\includegraphics[width=0.7\textwidth]{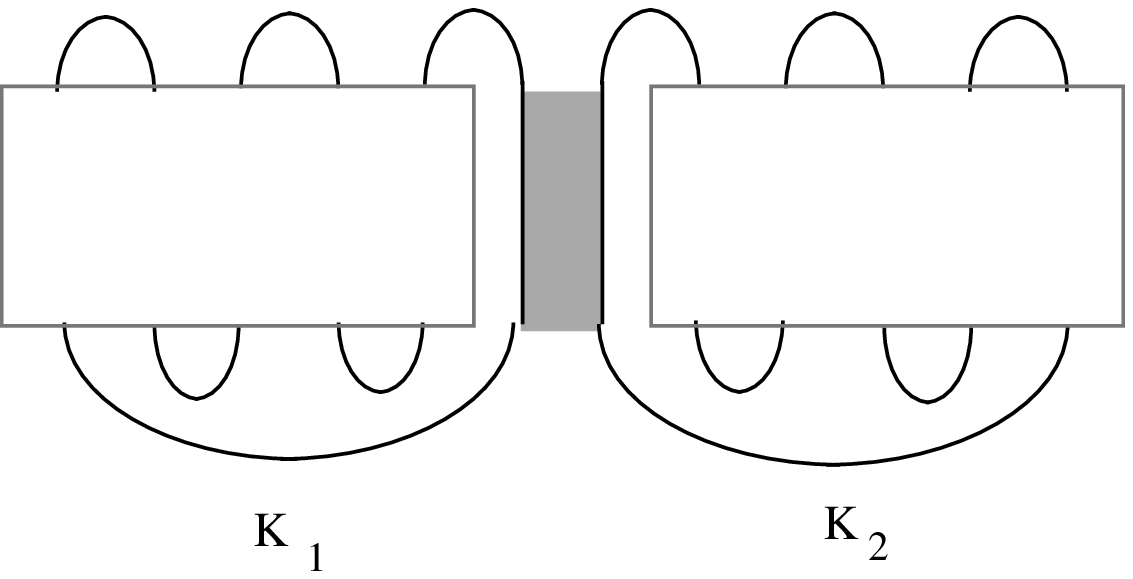}
\caption{} \label{fig:bridgeadd}
\end{figure}

To conclude this section, consider again the property of bridge number 
that first attracted Schubert:  its good behavior under summation of 
knots.  If $K_{1}, K_{2}$ are knots in minimal bridge 
position, there is a natural way to get a bridge-positioning of the 
sum $K_{1} \# K_{2}$.  Namely, arrange (as one can) that in the bridge 
positioning of the $K_{i}$, the right-most vertical strand of $K_{1}$ has 
no crossings and the left-most vertical strand of $K_{2}$ has no 
crossings.  Put $K_{1}$ to the left of $K_{2}$ and do the connected 
sum along these vertical strands.  The result is a bridge-positioning 
of $K_{1} \# K_{2}$ with $\bbb(K_{1}) + \bbb(K_{2}) - 1$ bridges; 
Schubert's theorem says that this is a minimal bridge positioning.  
See Figure \ref{fig:bridgeadd}.

\begin{figure} [tbh]
\centering
\includegraphics[width=0.5\textwidth]{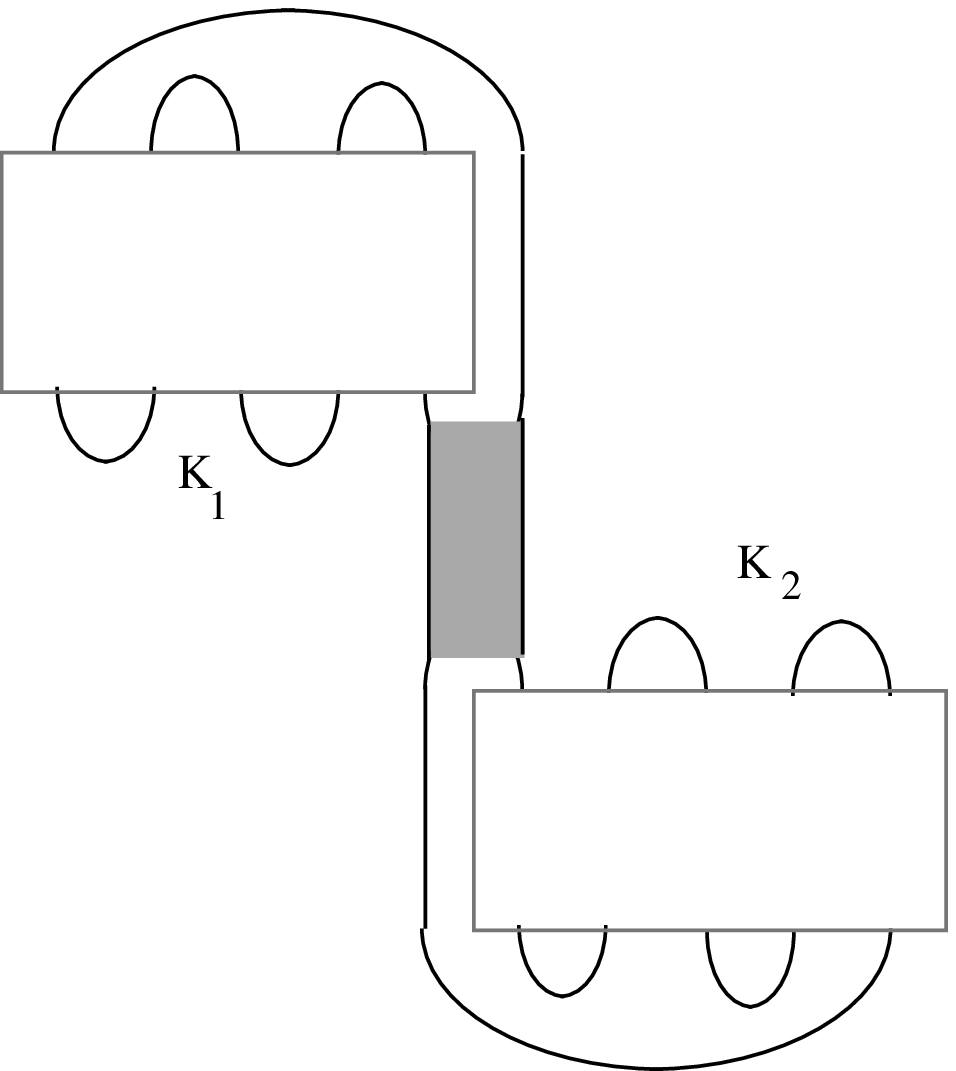}
\caption{} \label{fig:widthadd}
\end{figure}

There is a similar construction which one might easily conjecture
would do the same for width.  Put $K_{1}$ and $K_{2}$ in thin position
and position them so that $K_{1}$ lies entirely above $K_{2}$.  Then
do the connected sum of the knots via a monotonic band from the lowest
minimum of $K_{1}$ to the highest maximum of $K_{2}$.  See Figure
\ref{fig:widthadd}.  Of course we do not immediately know that the
result is a minimal width presentation for $K_{1} \# K_{2}$, but this
simple picture does show:

\begin{lemma} \label{lemma:widthadd}
    
    $$w(K_{1} \# K_{2}) \leq w(K_{1}) + w(K_{2}) - 2$$
    
    \end{lemma}
    
With the precedent of Schubert's theorem before us, it's natural to
ask
    
    \begin{quest} Is the inequality in Lemma \ref{lemma:widthadd} ever
    strict, or is it always true that $$w(K_{1} \# K_{2}) = w(K_{1}) +
    w(K_{2}) - 2?$$  In other words, does the positioning shown in 
    Figure \ref{fig:widthadd} always 
    minimize width?
    
    \end{quest}
    
Rieck and Sedgwick \cite{RS} show that the answer is yes when neither
$K_{1}$ nor $K_{2}$ has in its exterior an essential meridional planar
surface; Wu notes that their result now follows easily from Theorem
\ref{theorem:wu}.  See also and \cite{He}.  Without the assumption on
essential planar surfaces in the knot exteriors, only a little is
currently known, e.g. $$w(K_{1} \# K_{2}) \geq max\{ w(K_{1}),
w(K_{2}) \},$$ cf \cite{ScSc}, though there is strong evidence that
the answer is no, cf \cite{ST3}.

\section{From knots to graphs}

\subsection{Normal form for trivalent graphs}

Thin position techniques outlined above for knots $K \subset S^{3}$ do
not particularly make use of the fact that $K$ is connected (though
some of the applications do) and so thin position can be applied to
links as well as knots in $S^{3}$.  On the other hand, some thought is
needed if thin position is to be applied to imbedded graphs in
$S^{3}$.  We will restrict our discussion to trivalent graphs;
presumably more general graphs can be treated similarly but so far
there seems to be no notable application to higher valence graphs.

\begin{figure} [tbh]
\centering
\includegraphics[width=0.7\textwidth]{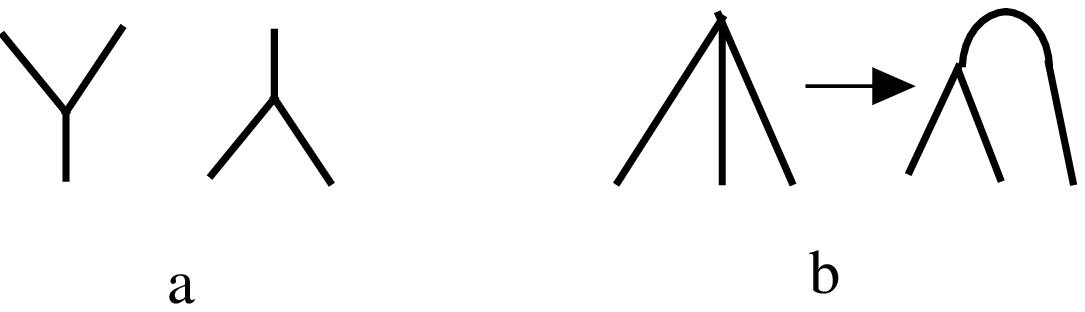}
\caption{} \label{fig:threemax}
\end{figure}

\begin{defin} \label{defin:normal} 
    {\rm Let $\Ggg$ be a finite trivalent
graph in $S^{3} - \{ poles \}$ and let $h: S^{3}
\map \mathbb{R}$ be the standard height function.  $\Ggg$ is in {\em normal
form} with respect to $h$ if

\begin{enumerate}
    
\item the critical points of $h|edges$ are nondegenerate and each lies 
in the interior of an edge; 

\item the critical points of $h|edges$ and the vertices of $\Ggg$ all
occur at different heights and

\item At each (trivalent) vertex $v$ of $\Ggg$ either two ends of 
incident edges lie above $v$ (we say $v$ is a {\em $Y$-vertex}) or two 
ends of incident edges lie below $v$ (we say $v$ is a {\em $\lll$-vertex}).
(See Figure \ref{fig:threemax} a.)

\end{enumerate} }
    \end{defin}
    
Any $\Ggg \subset S^{3}$ can be perturbed by a small isotopy to be
normal; for example, note that if three edges are incident to the same
vertex from below, then a small isotopy moves the end of one edge so
that is incident from above and has a maximum near the vertex.  That
is, such a vertex is replaced by an interior maximum adjacent to a
$\lll$-vertex.  (See Figure \ref{fig:threemax} b.)

Suppose $\Ggg$ is in normal form with respect to $h$. 

\begin{defin} 
{\rm The {\em maxima} of $\Ggg$ consist of all local maxima of $h|edges$
and all $\lambda$-vertices.  The {\em minima} of $\Ggg$
consist of all local minima of $h|edges$ and all $Y$-vertices.  

A maximum (resp.  minimum) that is not a $\lambda$-vertex (resp. 
$Y$-vertex) will be called a {\em regular} maximum (resp.  minimum). 
The set of all maxima and minima (hence including all the vertices of
$\Ggg$) is called the set of {\em critical points} of $\Ggg$.  The
heights of the critical points are called the {\em critical values} or
{\em critical heights}.

$\Ggg$ is in {\it bridge position} if there is a level sphere, called 
a {\em dividing sphere}, that lies above all 
minima of $\Ggg$ and below all maxima.  }
\end{defin}

\bigskip

There are two sorts of complications introduced when valence $3$
vertices are allowed: there is some subtlety in finding an appropriate
calculation of width; and surfaces that are properly imbedded in the
graph complement may behave in a less orderly fashion.  We treat each 
of these in turn:

\subsection{Width for graphs}

A naive way to define width is to proceed just as in the case of knots
or links: For each generic $-1 < t < 1$ let $w(t) = |\Ggg \cap
P_{t}|$, a function that increases by two at a regular minimum of $h$,
and by one at a $Y$-vertex.  Similarly $w(t)$ decreases by two at a
regular maximum of $h$, and by one at a $\lll$-vertex.  Now pick
generic heights $r_{1}, \ldots, r_{n}$, each between a distinct
adjacent pair of critical heights, and calculate the sum $\Ss_{i}
w(r_{i})$.  A complication that this definition introduces is this:
there are now two types of maxima in $\Ggg$, regular maxima and
$\lll$-vertices.  If one is pushed down past the other, the width
changes.  In particular, if a maximum of unknown type is pushed down
past other maxima and other minima, we do not {\it a priori} know that
the width is decreased.  (Symmetric statements are true, of course,
for minima.)  Such a definition then would make arguments using upper
and lower disks, arguments that worked so successfully in the case of
knots, pretty useless for graphs.

\begin{figure}[tbh]
\centering
\includegraphics[width=.9\textwidth]{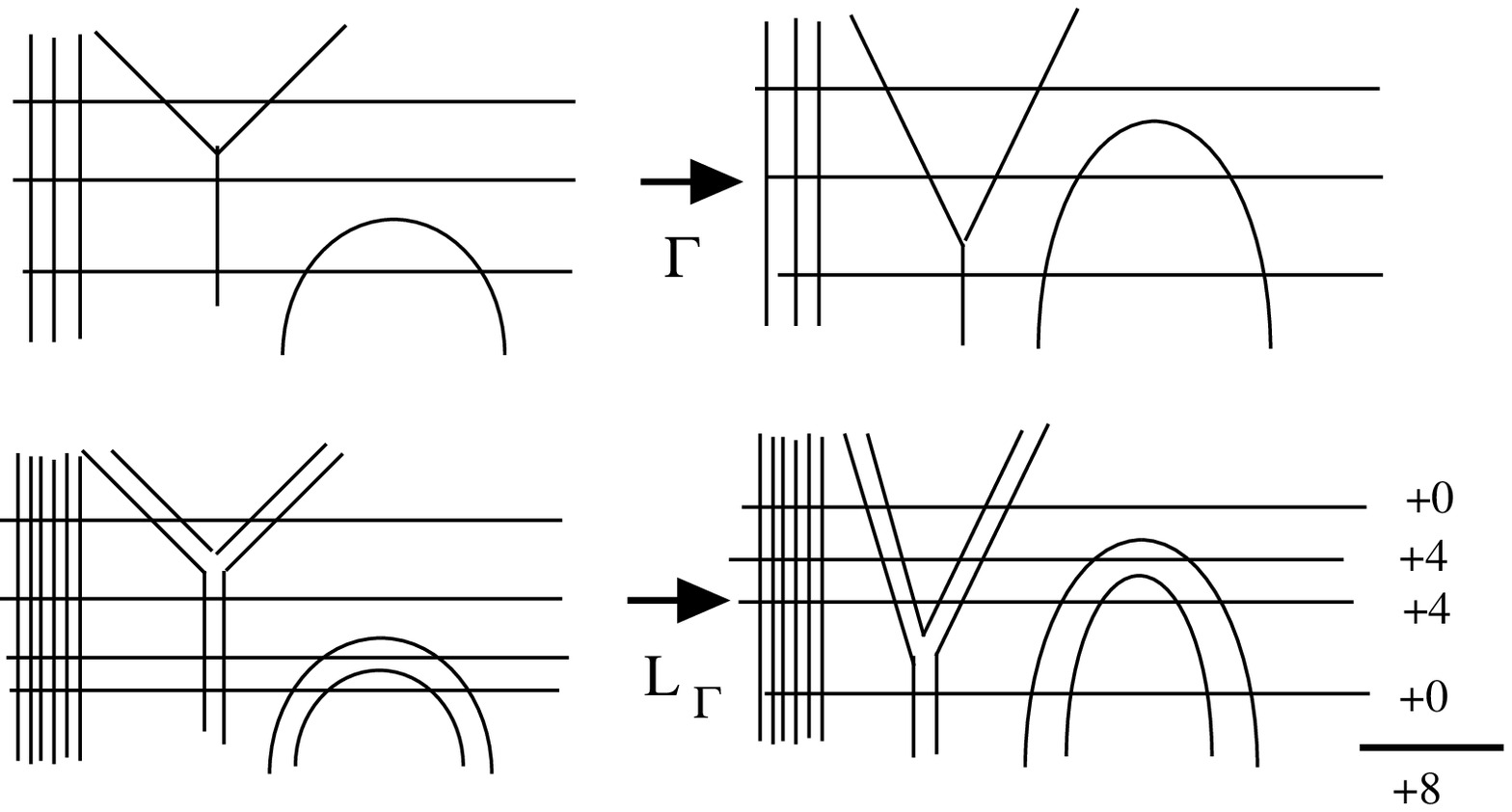}
\caption{} \label{fig:minmax1}
\end{figure}

A fix for this is to alter the definition slightly, taking into
account what the nature of the adjacent critical height is.  A
motivating thought is this: If we widen the graph to look like a
ribbon whose core is the graph, then the boundary of the ribbon is
just a standard link $L_{\Ggg}$, with no vertices.  Exchanging the
heights of two critical points of $\Ggg$ will rearrange heights of
several critical points on $L_{\Ggg}$ but in a predictable way. 
Moving maxima past maxima in $\Ggg$ will rearrange only maxima of
$L_{\Ggg}$ (and so have no effect on height); moving a maximum past a
minimum in $\Ggg$ will similarly move maxima past minima in
$L_{\Ggg}$.  See Figure \ref{fig:minmax1}.  It's relatively easy to
express the width of $L_{\Ggg}$ via the number of intersections of
$\Ggg$ with level planes, and we will use that as the definition of 
the width of $\Ggg$.  The upshot is a somewhat more complicated 
definition, but one that automatically has the property we seek: it's 
indifferent to the exchange in heights of adjacent maxima or adjacent 
minima, but will go down if a maximum is moved below a minimum.

\begin{defin} \label{defin:widthgraph} 
    {\rm Let $c_0 < \ldots < c_n$ be the
successive critical heights of $\Ggg$.  Let $r_i, 1 \leq i \leq n$ be
generic levels chosen so that $c_{i-1} < r_i < c_i$ and let $p_{i}$ 
denote the critical point at height $c_{i}$. For each $i, 1 
\leq i \leq n$ define $\rrr_{i}$ by:

\[\rrr_{i} = \left\{\begin{array}{cl}
                      2 & \mbox{when $p_{i-1}, p_{i}$ are both vertices}\\
		      3 & \mbox{when exactly one of $p_{i-1}, p_{i}$ is a vertex}\\
		      4 & \mbox{when $p_{i-1}, p_{i}$ are both regular 
		      critical points}
\end{array} \right. \]

Define the width of $\Ggg$ with respect $h$ to be $$W(\Ggg, h) =
\Sigma_{i} \rrr_{i} \cdot w_{r_i}.$$}
\end{defin}

This is a essentially the width introduced in \cite[Section 
3]{GST}, but made symmetric with respect to reflection through a 
horizontal plane.

\begin{defin} 
{\rm A graph $\Ggg \subset S^3$ is in thin position (with respect to $h$)
if $W(\Ggg, h)$ cannot be lowered by an isotopy of $\Ggg$.  In that 
case, $W(\Ggg, h)$ is denoted $W(\Ggg)$.}
\end{defin}

\begin{figure}[tbh]
\centering
\includegraphics[width=.9\textwidth]{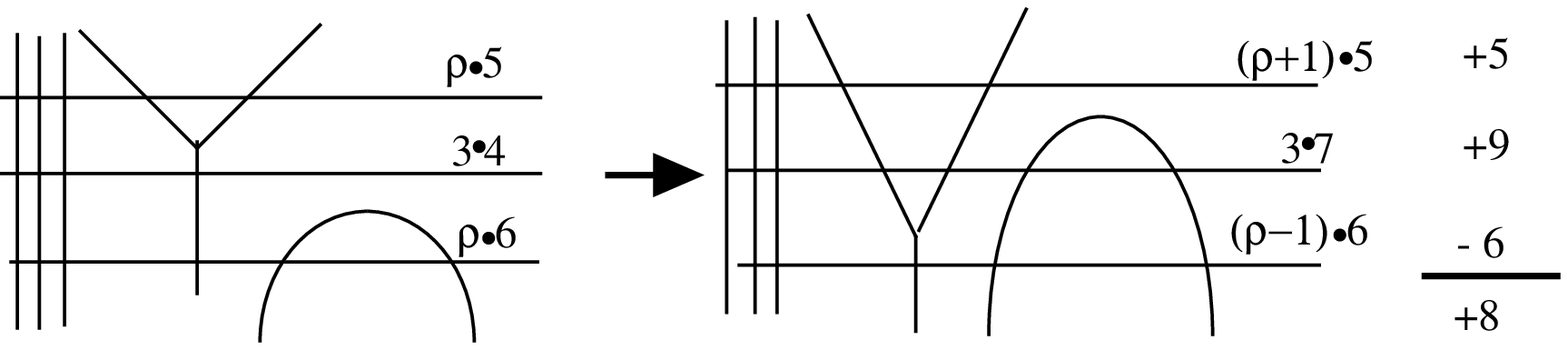}
\caption{} \label{fig:minmax2}
\end{figure}

{\bf Remark:} In practice, the chief property of the width $W(\Ggg,
h)$ that we will need is this: The width is decreased if a maximum is
pushed below a minimum, but the width is unaffected by pushing one
maximum above or below another maximum, or one minimum above or below
another minimum.  See Figures \ref{fig:minmax1}, \ref{fig:minmax2}.

\bigskip

\subsection{Surfaces in graph complements}

Surfaces that lie in a knot complement have the pleasant feature that 
boundary components are either horizontal (if the boundary component 
is a meridian) or the height function on the boundary circles roughly 
follow that of the knot, having maxima where the knot has a maximum 
and similarly with minima.  It's easy to locally isotope the surface 
so that all singularities of $h$ on its boundary are local 
half-centers.

\begin{figure}[tbh]
\centering
\includegraphics[width=.9\textwidth]{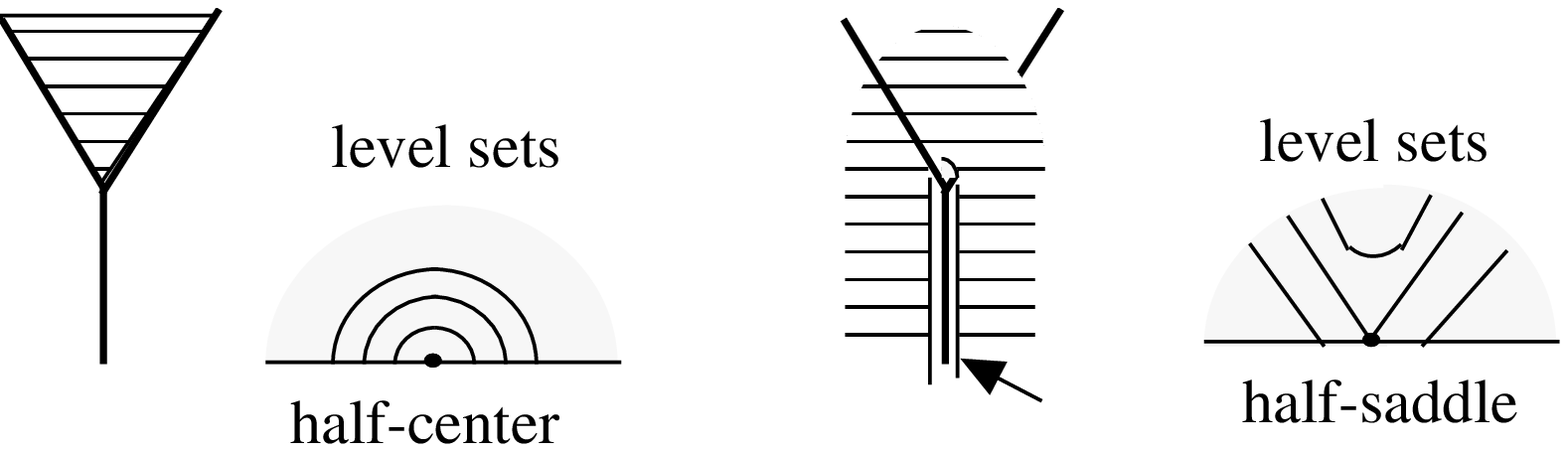}
\caption{} \label{fig:saddlecenter}
\end{figure}

In contrast, surfaces in a graph complement may have minima at
$\lll$-vertex maxima and maxima at $Y$-vertex minima.  See Figure
\ref{fig:saddlecenter}.

A graph $\Ggg \subset S^3$ in normal form with respect to $h$ can be
thickened slightly to give a solid handlebody $\eta(\Ggg) \subset S^3$
with the predictable height structure (e.g.  very near any regular
maximum of $\Ggg$ there are two non-degenerate critical points of
$h|\bdd(\eta(\Ggg))$, one a saddle just below and one a maximum just
above.)  We will be concerned with simple closed curves on
$\bdd\eta(\Ggg)$ and with properly imbedded surfaces in $S^3 -
\eta(\Ggg)$.

\begin{defin} {\rm Suppose $\Ggg$ is a graph, in normal form with respect
to $h$, and $c \subset \bdd\eta(\Ggg)$ is a simple closed curve.  Then
$c$ is in normal form on $\bdd\eta(\Ggg)$ if either it is a horizontal
meridian circle or each critical point of $h$ on $c$ is
non-degenerate, and occurs near an associated critical point of $\Ggg$
in $\bdd\eta(\Ggg)$.  Furthermore, the number of critical points of
$c$ has been minimized via isotopy of $c$ in $\bdd\eta(\Ggg)$. }
\end{defin}

\begin{defin} \label{def:normal} {\rm A properly imbedded surface 
$$(F,\bdd F) \subset (S^3 - \eta(\Ggg), \bdd\eta(\Ggg))$$ is in {\em 
normal form} if
\begin{enumerate}

\item each critical point of $h$ on $F$ is nondegenerate, 

\item $\bdd F$ is in normal form with respect to $h$

\item no critical point of $h$ on int$(F)$ occurs near a 
critical height of $h$ on $\Ggg$, 

\item no two critical points of $h$ on int$(F)$ or $\bdd F$ occur at
the same height,

\item the minima (resp.  maxima) 
of $h|\partial F$ at the minima (resp.  maxima) of $\Ggg$ are 
also local extrema of $h$ on $F$, i.e., `half-center' singularities, 

\item the maxima of $h|\partial F$ at $Y$-vertices and the minima of 
$h|\partial F$ at $\lambda$-vertices are, on the contrary, 
`half-saddle' singularities of $h$ on $F$.
\end{enumerate} }
\end{defin}

Standard Morse theory ensures that, for $\Ggg$ in normal form, any
properly imbedded surface $(F, \partial F)$ can be properly isotoped
to be in normal form.  

\bigskip

\subsection{Upper and lower triples}

The definition of upper and lower disks 
naturally extends to the context of graphs:

\begin{defin} \label{defin:graphupper}
{\rm Given $\Ggg$ in normal form and $P$ a level sphere for $h$ at a
generic height, let $B_u$ and $B_l$ denote the balls which are the
closures of the region above $P$ and below $P$ respectively.  An {\em
upper disk} (resp.  {\em lower disk}) for $P$ is a disk $D \subset S^3
- \eta(\Ggg)$ transverse to $P$ such that $\partial D =
\alpha\cup\beta$, where $\beta$ is a normal arc imbedded on
$\partial\eta(\Ggg)$, $\alpha$ is an arc properly
imbedded in $P - \eta(\Ggg)$, $\partial\alpha=\partial\beta$, and a
small product neighborhood of $\aaa$ in $D$ lies in $B_u$ (resp. 
$B_l$) i.e., it lies {\em above} (resp.  {\em below}) $P$. }
\end{defin}

Note that $D \cap P$ consists of simple closed curves and arcs with
ends in $\bbb$.  A natural occurence of upper (or, symmetrically,
lower) disks is this: According to Definition \ref{def:normal}, a
maximum of $\bdd F$ near a maximum of $\Ggg$ is a half-center
singularity on $\bdd F$.  In particular, a sphere $P$ just below this
maximum will cut off an upper disk from $F$. 

As was noted above, curves in $\partial\eta(\Ggg)$ can be quite
complicated, so a somewhat more elaborate notion than upper or lower
disk will be needed.

\begin{defin} \label{defin:triple} {\rm Suppose, as above, $\Ggg$ is in
normal form and $P$ is a level sphere for $h$ at a generic height.  An
{\em upper triple} (resp.  lower triple) $(v, \aaa, E)$ for $P$ is an upper (resp. 
lower) disk $E$ with these properties

\begin{enumerate}

    \item The arc $\aaa \subset \bdd E$ of Definition
    \ref{defin:graphupper} has its ends at different points
    of $\Ggg \cap P$ (i.e. $\aaa$ is not a loop)
     
     \item $v$ is one of the points of $\Ggg \cap P$ at an end of 
     $\aaa$, and
     
     \item although there may be arc components of $\inter(E) \cap P$,
     none of them is incident to $v$.

\end{enumerate} }
\end{defin}

For example, in the old setting of, say, $F$ a Seifert surface for a
knot $K$, any arc $\aaa$ of $F \cap P$ that is inessential in $F$ cuts
off either an upper or a lower disk $E$.  If $v$ is either point of $K
\cap P$ at the ends of $\aaa$, then $(v, \aaa, E)$ is an upper or
lower triple.

\subsection{An application}

As an illustration of how thin position can be used for graphs -- in
particular, why it is useful to have a definition that is indifferent
to pushing maxima past maxima -- we'll offer an updated proof of the
key Theorem in \cite{ST1}, which leads to the proof that any Heegaard
splitting of $S^{3}$ is standard.  (The roots of this proof go back to
Otal \cite{Ot}.)  For the proof in \cite{ST1} we did not have in hand
the efficient Definition \ref{defin:widthgraph} of width; instead we
used a rather clumsy alternative, examining the entire function $w(t)$
and minimizing its maximum, together with the number of times the
function achieves that maximum.

The setting we consider is this: $\Ggg \subset S^{3}$ is a finite
graph in normal form whose complement is $\bdd$-reducible, so there is
a disk $(D, \bdd D) \subset (S^{3} - \eta(\Ggg), \bdd \eta(\Ggg))$ in
which $\bdd D$ is essential in $\bdd \eta(\Ggg)$.  The edges of the
graph $\Ggg$ are allowed to slide over each other.

\begin{lemma} \label{lemma:isatriple} Suppose $P$ is a generic level
sphere for $\Ggg \subset S^{3}$, $S^{3} - \eta(\Ggg)$ is
$\bdd$-reducible, and the $\bdd$-reducing disk $D$ has been chosen to
minimize $|D \cap P|$ and is properly isotoped to be in normal form. 
Then either there is an edge of $\Ggg$ that is disjoint from $\bdd D$
or there is a point $v \in P \cap \Ggg$ with the following property:
Suppose $\aaa$ is an arc of $P \cap D$ that is outermost among the
set of arcs of $P \cap D$ that are incident to $v$.  Let $E$ be the
disk it cuts off from $D$.  Then $(v, \aaa, E)$ is either an upper or
a lower triple.
\end{lemma}

\begin{proof} Since $D$ was chosen to minimize $|D \cap P|$ it follows 
that every component of $D \cap P$ is essential in the planar 
surface $P - \eta(\Ggg)$.  If any intersection point of $\Ggg$ with 
$P$ is incident to no arc component of $D \cap P$ then the edge of 
$\Ggg$ containing that point is disjoint from $\bdd D$ and we are 
done.  So we may as well assume that each point of $\Ggg \cap P$ is 
incident to some arc component of $D \cap P$; it follows that some 
point $v \in \Ggg \cap P$ is incident to no loops at all.  Among all 
arcs of $D \cap P$ that are incident to $v$, let $\aaa$ be the arc 
that is outermost on $D$.  Then by construction, $\aaa$ cuts off an 
upper disk $E$ (say) from $D$ in which no other arc of intersection is 
incident to $v$ and $\aaa$ is not a loop in $P$.
\end{proof}

\begin{lemma} \label{lemma:graphbridge} Suppose the edges of $\Ggg$ 
have been slid and isotoped so as to minimize $W(\Ggg)$ (cf Definition 
\ref{defin:widthgraph}) and suppose $S^{3} - \eta(\Ggg)$ is 
$\bdd$-reducible.  Then either the edges of $\Ggg$ can be slid until 
there is a $\bdd$-reducing disk whose boundary is disjoint from an 
edge, or $\Ggg$ is in bridge position.
\end{lemma}

\begin{proof} Suppose on the contrary that $\Ggg$ is not in bridge
position.  Let $P$ be a thin level sphere, i.  e. a level sphere
intersecting $\Ggg$ so that the adjacent critical heights above and
below $P$ are a minimum (possibly a $Y$-vertex) and a maximum
(possibly a $\lll$-vertex) respectively.  Choose a $\bdd$-reducing
disk $D$ so as to minimize $|D \cap P|$.  If any edge of $\Ggg$ is not
incident to $\bdd D$ then we are done.  If every edge is incident to
$\bdd D$ then in particular $D \cap P$ is incident to every point in
$\Ggg \cap P$.  In that case, let $(v, \aaa, E)$ be the upper (say)
triple given by Lemma \ref{lemma:isatriple}.  Then $E$ may be used to
slide one end $\eee$ of the edge of $\Ggg$ on which $v$ lies down to
$\aaa \subset P$.  (The details of this move, involving possible
``broken edge slides'' are a bit more complicated than it might first
appear, cf.  \cite[Proposition 2.2]{ST1}.)

Unfortunately this move, so similar to the one used in Theorem
\ref{theorem:arcthin}, does not in this case necessarily thin $\Ggg$. 
To see why, suppose that, before the slide, the end $\eee$ simply
ascends from $v$ into a $Y$-vertex, from below.  Then the slide we've
just described will move the end to $\aaa$ (creating a $\lll$-vertex
just below $P$ and lowering $|P \cap \Ggg|$ by $1$); but also the
$Y$-vertex merely becomes a regular minimum.  Both the new
$\lll$-vertex and the transformation of the $Y$-vertex into a regular
minimum will actually raise the width, by a total of $4|P \cap \Ggg| -
2$.  See Figure \ref{fig:newvertex}.  This is a technical setback, but
not a devastating one, as we now briefly outline.

Note that the slide we've just described lowers $w_{P} = |P \cap
\Ggg|$, so we cannot repeat the process indefinitely.  The argument
stops either because there is a $\bdd$-reducing disk disjoint from
an edge (and we are done) or when $P$ is no longer a thin sphere.  In
the latter case, either all the minima above $P$ or all the maxima
below $P$ have been removed by the sequence of edge slides.  With no
loss, assume that the process stops because all the minima above $P$
have been removed.  (We do not assume that an upper disk arises at
each stage, but of course an upper disk does happen to be needed at
the last stage, since a lower disk would give rise to a minimum just
above $P$).  We will show that by the time the process stops, the
width has been reduced.

\begin{figure}[tbh]
\centering
\includegraphics[width=.9\textwidth]{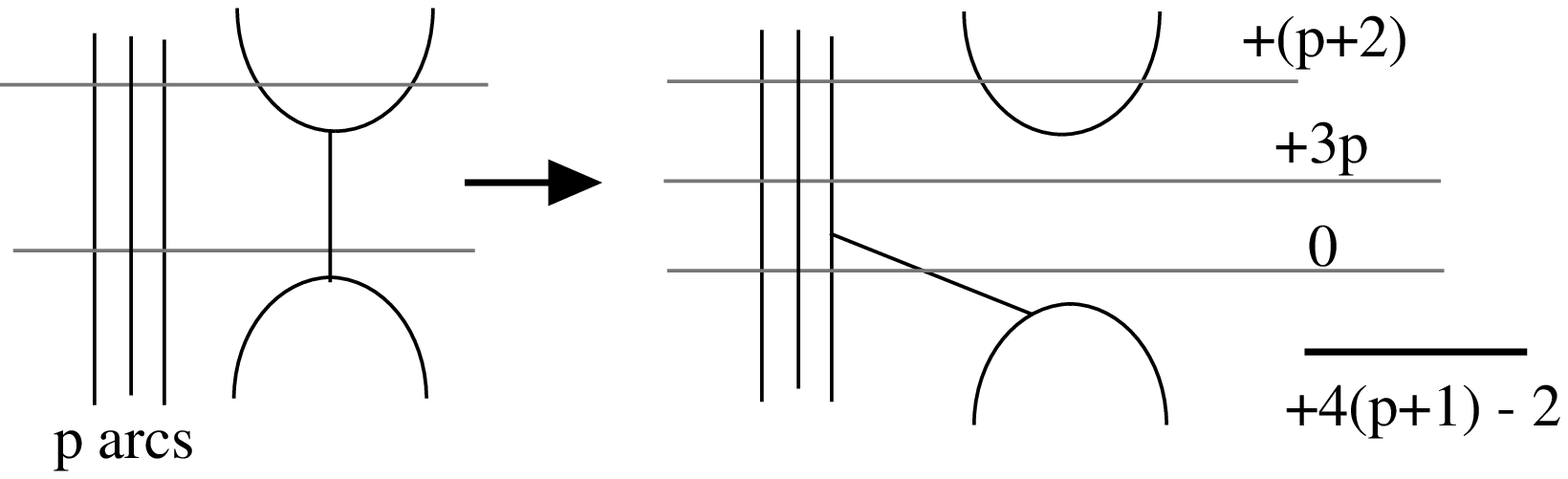}
\caption{} \label{fig:newvertex}
\end{figure}

We will assume that there are no further thin levels above $P$ and
leave it to the reader to adjust the following argument for the
general case (by counting for $a$ and $b$ only minima just above $P$
and then subtracting a further such term for each thin sphere above
$P$).  Let $a$ be the number of regular minima and $b$ be the number
of $Y$-vertices lying above $P$.  Define $$W_{P}(\Ggg) = W(\Ggg) - 4(2a
+ b)w_{P}.$$

{\bf Claim:} Every move in the process (whether on an upper or a lower
disk) decreases $W_{P}(\Ggg)$.

{\bf Proof of claim:} Each of the moves is either a simple isotopy of
an arc of $\Ggg$ to an arc lying just above or below $P$, or it is a
slide of an end $\eee$ of an edge of $\Ggg$ to an arc just above or
below $P$.  In the former case, no vertex moves and the proof is
almost immediate.  The move reduces $W_{\Ggg}$ and we only have to
check that the reduction is greater than the increase of $-4(2a +
b)w_{P}$ which may result from the elimination of regular minima of
$\Ggg$ that lie above $P$.  These regular minima necessarily lie on
the arc $\bbb \subset \Ggg$ being isotoped.  For every minimum of
$\bbb$ above $P$ there will be a maximum (and {\it in toto} one more maximum
than minimum on each component of $\bbb - P$ that lies above $P$).
Eliminating both a minimum and a maximum above $P$ will reduce
$W(\Ggg)$ by a total of at least $8w_{P} + 16$.  (Since width is not
altered by rearranging orders of maxima or of minima, we may assume
for the purposes of calculation that the cancelling critical points are
respectively the lowest maximum and the highest minimum, i.e. at
adjacent heights).  At the same time, $-4(2a + b)w_{P}$ will go up by
only $8w_{P}$.  Thus, in any case, the isotopy of $\bbb$ will reduce
$W_{P}(\Ggg)$.

Next suppose the move is an edge slide.  If the end $\eee$ of an edge
that is slid descends into a $Y$-vertex or ascends into a
$\lll$-vertex, then the edge slide does not create more critical
levels and again the argument is fairly straightforward: $W(\Ggg)$ is
always reduced, and if $\eee$ ascends into a $\lll$-vertex or if it
descends into a $Y$-vertex lying below $P$ then the only minima above
$P$ that disappear are internal minima for which the above argument
applies.  If the end of $\eee$ descends into a $Y$-vertex above $P$,
that minimum is eliminated (raising $-4(2a + b)w_{P}$ by $4w_{P}$) but
it can be viewed as being cancelled with an adjacent maximum, which
reduces $W(\Ggg)$ by at least $4w_{P} + 6$.

So, not surprisingly, to prove the claim we are reduced to the case in
which an extra critical level may be created, because the end $\eee$
either descends into a $\lll$-vertex or ascends into a $Y$-vertex. 
(In these cases the slide does not eliminate the critical point at its
end, and creates a new one near $P$.)  Let us call the terminating
vertex $v$.  There are four cases: $v$ may lie above or below $P$ and
the disk defining the move may be an upper or lower one.  But, for
example, if $v$ lies above $\eee$ and the disk is a lower one, we may
imagine the slide as the composition of one based on an upper disk
taking $v$ down below $P$, followed by one based on a lower disk
bringing $v$ back up to $P$.  In other words, it suffices to consider
the two cases where $v$ lies below (resp.  above) $P$ and the disk
determining the slide is a lower disk (resp.  upper disk).

Suppose first that $v$ is below $P$ and the move is via a lower disk,
so the vertex is moved to a $Y$-vertex just above $P$.  We have
already argued that eliminating internal critical points on $\eee$ can
only improve the situation, so we may as well assume that $\eee$
either descends from $P$ straight down into a $\lll$-vertex or $\eee$
has a single internal minimum, adjacent to its ascent from below into
a $Y$-vertex.  In the latter case the internal minimum of $\eee$ is
also eliminated so again no new critical level is really created: In
fact the slide is equivalent to a move that just brings the $Y$ vertex
up to $P$.  Since raising a $Y$-vertex cannot raise the width, and a
new $Y$-vertex above $P$ reduces $-4(2a+b)$, the overall effect is to
reduce $W_{P}(\Ggg)$.  In the case where the end $\eee$ descends
straight into $v$, a $\lll$-vertex, the slide that moves $v$ up to $P$
can be viewed as the composition of two moves: first move $v$ up to a
new $Y$-vertex just above the lowest thin sphere $P'$ above $v$.  We
have already seen (cf.  Figure \ref{fig:newvertex}) that this raises
the width by $4w_{P'} - 2$.  Next raise $v$ up above $P$.  This lowers
the width every time $v$ passes a maximum: by $4$ for every
$\lll$-vertex passed and by $8$ for every regular maximum passed (cf. 
Figure \ref{fig:minmax2}).  But then the total amount of the
reduction, determined by the number and type of maxima between $P'$
and $P$ is at least $4(w_{P'} - w_{P})$ so in the end the move raises
the width by at most $4w_{P} - 2$.  When combined with the effect of
raising $b$ by $1$, the result is that $W_{P}(\Ggg)$ goes down by at
least $2$ (indeed exactly $2$ only if $P' = P$).

Finally, suppose $v$ is above $P$ and the move is via an upper disk. 
Again we may as well assume, because eliminating internal critical
points only improve the situation, that $\eee$ either simply ascends
from $P$ into a $Y$-vertex or its terminating end descends into a
$\lll$-vertex from a single adjacent internal maximum.  We have seen
that in the former case the width increases by exactly $4w_{P} - 2$. 
On the other hand, the $Y$-vertex at $v$ becomes a regular vertex,
reducing $-4(2a + b)$ by $4$.  The net effect is to reduce
$W_{P}(\Ggg)$ by $2$.  In the latter case, again the slide effectively
just moves the $\lll$-vertex below $P$, reducing the width by moving a
$\lll$-vertex down past minima, but having no effect on $a$ or $b$. 
This finally proves the claim in all cases.

Now let $\Ggg'$ be the graph when the process stops, with no further
minima above $P$.  We have just seen $W_{P}(\Ggg') < W_{P}(\Ggg) <
W(\Ggg)$.  But since there are no minima of $\Ggg'$ above $P$,
$W_{P}(\Ggg') = W(\Ggg')$.  Hence $W_{P}(\Ggg') < W(\Ggg)$, a
contradiction to the original assumption that $\Ggg$ was in thin 
position.
\end{proof}

\begin{lemma} \label{lemma:aretwotriples} Suppose $\Ggg$ is in bridge 
position and $S^{3} - \eta(\Ggg)$ is $\bdd$-reducible.  Then the edges 
of $\Ggg$ can be slid rel a dividing sphere until either

\begin{itemize}

\item there is a $\bdd$-reducing disk whose boundary is disjoint from an edge or

\item for some dividing sphere $P$ there are both upper and lower 
triples $(v_{u}, \aaa_{u}, E_{u})$, $(v_{l}, \aaa_{l}, E_{l})$ so that 
the disks $E_{u}, E_{l}$ are disjoint in $S^{3} - \eta(\Ggg)$ 
\end{itemize}

In the latter case, the triples may further be chosen so that either
no arc of $(E_{u} - \aaa_{u}) \cap P$ is incident to $v_{l}$ or, vice
versa, no arc of $(E_{l} - \aaa_{l}) \cap P$ is incident to $v_{u}$.
\end{lemma}

\begin{proof} Let $P$ be a dividing sphere and choose $D$ among all 
$\bdd$-reducing disks for $S^{3} - \Ggg$ so that $|\bdd D \cap P|$ is 
minimal.  In particular, this guarantees that no arc component of $D 
\cap P$ is a trivial loop in the planar surface $P - \Ggg$.

\bigskip

{\bf Claim:} $\Ggg$ can be slid and isotoped rel $P$ so that it is
still in bridge position and at the lowest maximum (resp.  the highest
minimum) of $\Ggg$, $\bdd D$ also has a maximum (resp.  minimum).  So
for $D$ in normal form, both the lowest maximum and highest minimum
of $\Ggg$ are incident to half-center singularities on $\bdd D$.

\bigskip

{\bf Proof of Claim:} Choose any component $\Ggg_{0}$ of $\Ggg - P$,
say one lying above $P$.  Since $\Ggg$ is in bridge position,
$\Ggg_{0}$ is necessarily a tree.  In particular, $\eta(\Ggg_{0}) - P$
is a planar surface.  If $\bdd D$ is not incident to an edge of
$\Ggg_{0}$ we are done, so assume it is incident to every edge; it
follows (by examining an innermost loop, if any, on the planar surface
$\eta(\Ggg_{0}) - P$) that there is a component of $\bdd D - P$ that
runs from the end $\eee_{1}$ of one edge of $\Ggg_{0}$ to the end of
another $\eee_{2}$.  Imagine collapsing the edges of $\Ggg_{0}$ that
are not incident to $P$ to a single vertex (so $\Ggg_{0}$ is simply
the cone on its ends) then sliding to recreate a trivalent graph with
only maxima, in which a single pair of edges (forming a maximum that
we may isotope to be the lowest maximum) contains the entire subarc of
$\bdd D$ that connects $\eee_{1}$ to $\eee_{2}$.  See Figure
\ref{fig:bridgegraph}.  The new graph is again in bridge position and
is homeomorphic to the original, so the width has not been altered. 
This establishes the claim.

\begin{figure}[tbh]
\centering
\includegraphics[width=.9\textwidth]{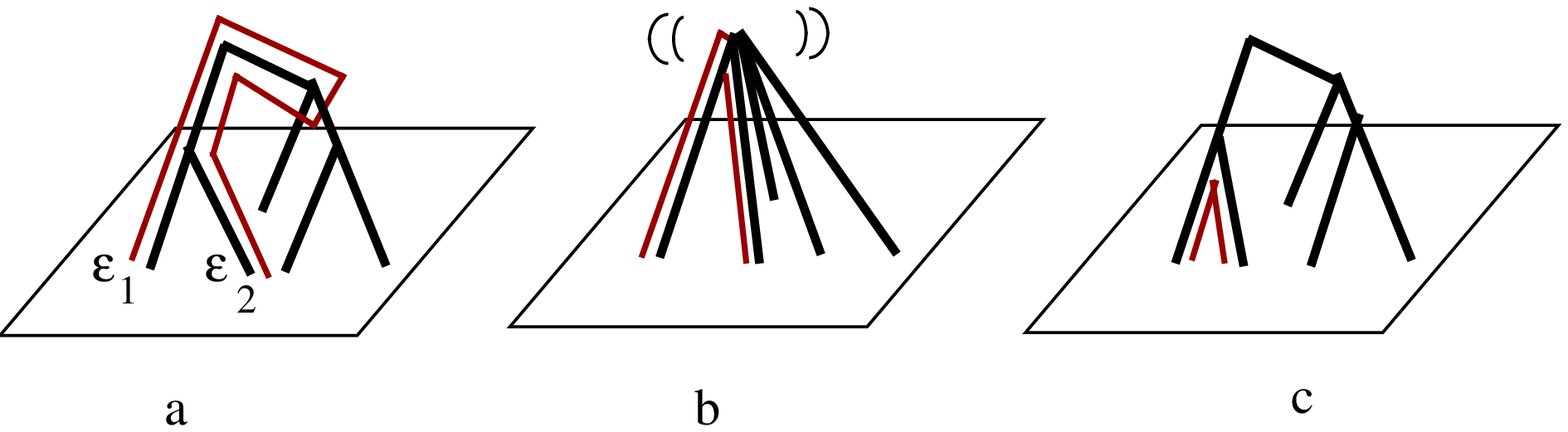}
\caption{} \label{fig:bridgegraph}
\end{figure}

Following the Claim, note that a level sphere just below the lowest
maximum will cut off an upper disk entirely contained above the
sphere; moreover the arc of intersection with $P$ is not a loop. 
Hence it's an upper triple.  Similarly, a level sphere just above the
highest minimum will cut off a lower triple.  According to Lemma
\ref{lemma:isatriple} every generic level sphere in between cuts off
either an upper or a lower triple.  So, as usual, there is a level
sphere that cuts off both an upper and a lower triple, $(v_{u},
\aaa_{u}, E_{u})$ and $(v_{l}, \aaa_{l}, E_{l})$.  Moreover at least
one of the two, say $(v_{u}, \aaa_{u}, E_{u})$, is obtained via Lemma
\ref{lemma:isatriple}.  If, among the arcs incident to $v_{u}$, there
is also an arc cutting off a lower triple, use this triple for $(v_{l}
= v_{u}, \aaa_{l}, E_{l})$.  Then automatically no arc of $(E_{l} -
\aaa_{l})$ is incident to $v_{u}$, establishing the last property
required.  If, on the other hand, no arc incident to $v_{u}$ cuts off
a lower triple, then every arc outermost among the arcs incident to
$v_{u}$ must cut off an upper triple.  In this case, to establish the
last property of the lemma, suppose on the contrary that some arc of
$(E_{l} - \aaa_{l})$ is incident to $v_{u}$. Then an outermost such arc
will cut off a (possibly different) upper triple $(v'_{u}, \aaa'_{u},
E'_{u})$ with the property that no arc of $(E'_{u} - \aaa'_{u}) \cap
P$ is incident to $v_{l}$.
\end{proof}

\begin{lemma} Suppose $\Ggg$ is not the unknot, and the edges of 
$\Ggg$ have been slid and isotoped so as to minimize $W(\Ggg)$.  
Suppose further that $S^{3} - \eta(\Ggg)$ is $\bdd$-reducible.  Then 
the edges of $\Ggg$ can be further slid until there is a 
$\bdd$-reducing disk whose boundary is disjoint from an edge.
\end{lemma}

\begin{figure}
\centering
\includegraphics[width=.9\textwidth]{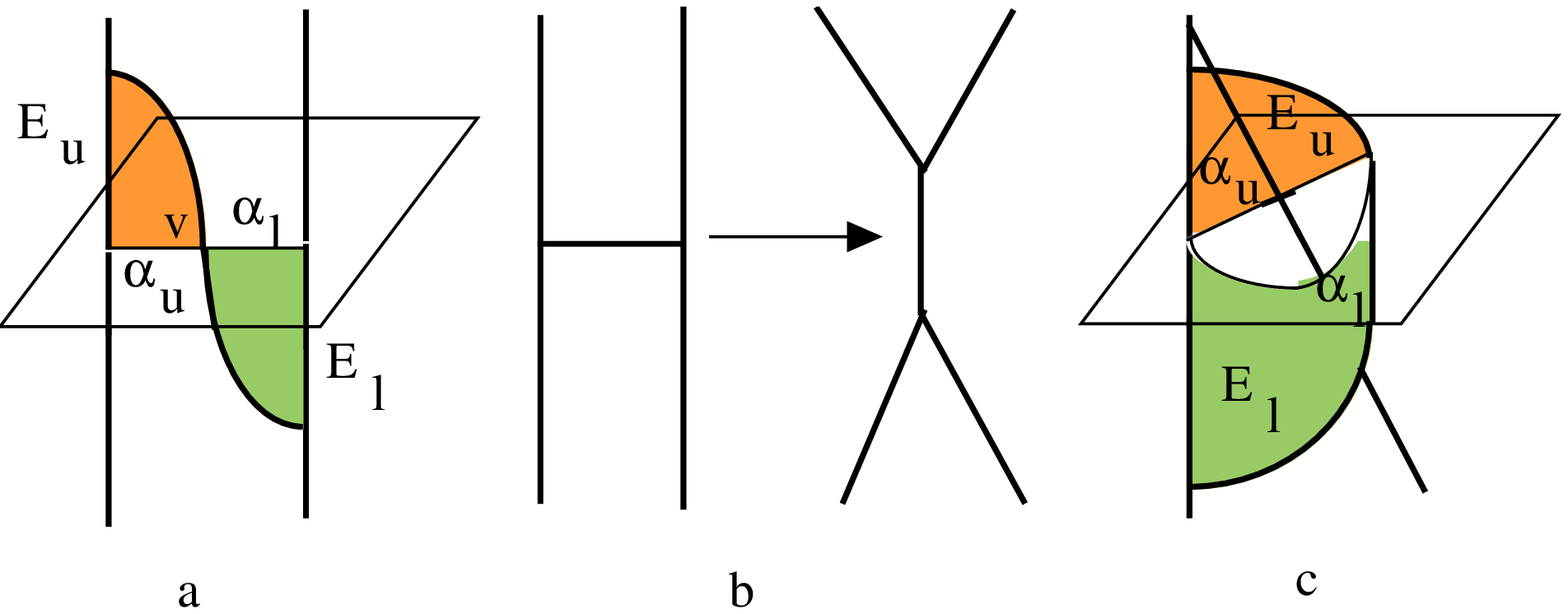}
\caption{} \label{fig:thintwo}
\end{figure}

\begin{proof} Suppose not.  Following Lemma \ref{lemma:graphbridge},
we may assume $\Ggg$ is in bridge position.  Then consider the upper
and lower triples $(v_{u}, \aaa_{u}, E_{u})$, $(v_{l}, \aaa_{l},
E_{l})$ given by Lemma \ref{lemma:aretwotriples} with respect to a
dividing sphere $P$.  In particular, we assume with no loss that no
arc of $(E_{u} - \aaa_{u}) \cap P$ is incident to $v_{l}$.  Then
$E_{u}$ may be used to slide an end of the edge on which $v_{u}$ lies
down to $\aaa_{u}$ without affecting the end of the edge on which
$v_{l}$ lies, so afterwards the latter end can also be brought to $P$. 
But sliding one end down and the other end up will typically reduce
the width, which is impossible.  An alternate possibility is that the
two moves actually level an entire edge, but again this would allow
the graph to be thinned.  See Figure \ref{fig:thintwo} a), b).  The
final possibility is that the two slides together simultaneously level
two edges (when $\aaa_{l}$ and $\aaa_{u}$ have the same pair of end
vertices), moving a cycle $\ggg \subset \Ggg$ in $\Ggg$ onto $P$. 
See Figure \ref{fig:thintwo} c).  

If either of the disk components of $P - \ggg$ is disjoint from
$\Ggg$, then either $\Ggg = \ggg$ (and we are done) or that disk
component is a $\bdd$-reducing disk as required.  But even if a disk
component $P_{0}$ of $P - \ggg$ intersects $\Ggg$, we can just apply
to $P_{0}$ the process we earlier applied to all of $P$ to find a
series of edge slides that will either finally exhibit an edge
disjoint from a $\bdd$-reducing disk $D$ (via a point in $\Ggg \cap
P_{0}$ incident to no arc of $P_{0} \cap D$) or will iteratively
reduce the number of points in $\Ggg \cap P_{0}$ until $\Ggg \cap
P_{0} = \emptyset$ so, as above, $P_{0}$ is a $\bdd$-reducing disk
disjoint from an edge.
\end{proof}

\begin{cor}  \label{cor:standard} Any Heegaard splitting of $S^{3}$ is 
standard.
\end{cor}

\begin{proof} Given a Heegaard splitting of $S^{3}$, let $\Ggg$ be a
spine of one of the handlebodies.  Apply the above argument not to
just a single $\bdd$-reducing disk for the complement (i.  e. a single
meridian disk for the complementary handlebody) but to a complete
collection of such disks.  The argument is essentially the same and
terminates either with an edge $e$ disjoint from a complete collection
or with $\Ggg$ simply the unknot, i.e. the standard genus one
splitting of $S^{3}$.  In the former case, a meridian circle of the
edge $e$ is disjoint from a complete collection of meridians for the
handlebody $S^{3} - \eta(\Ggg)$, so it also bounds a disk in $S^{3} -
\eta(\Ggg)$.  Thus there is a sphere in $S^{3}$ intersecting $\Ggg$ in
a single point in $e$.  This is a reducing sphere for the Heegaard
splitting which divides the splitting into two separate splittings of
the $3$-sphere.  The conclusion follows by induction on the genus of
the Heegaard splitting.
\end{proof}

\section{From graphs back to knots}

One might hope that thin position would be helpful in understanding
the tunnel structure for knots in $S^{3}$.  We noted above that for a
tunnel number one knot, thin position is bridge position.  That is, 
if $K$ has tunnel number one and is in thin position, then there are no
thin level spheres for $K$.  It's natural to ask about the behavior of the
tunnel arc with respect to the standard height function, once the knot
is in thin position.  The union of the knot $K$ and the tunnel arc
$\ttt$ is of course a graph in $S^{3}$; moreover, one way of viewing
the definition of unknotting tunnel is that the graph $\Ggg = K \cup
\ttt$ is the spine of a genus two Heegaard splitting of $S^{3}$.  The 
reason that this does not just fall into the program leading to 
Corollary \ref{cor:standard} is that in the knot tunnel case, whereas 
we are allowed to slide the ends of the tunnel over the knot, and 
over the other end of the tunnel, we can never regard a subarc of $K$ 
as an edge of $\Ggg$ that can be slid.  Nonetheless, the answer is 
simple and direct: the tunnel may be made level with respect to the 
standard height function and this is the thinnest positioning of 
$\ttt$ possible. 

\begin{theorem} \label{theorem:thin=level} Suppose that $K$ is a knot
with unknotting tunnel $\ttt$ and $K$ is in thin position with respect
to the standard height function $h$.  Then $\ttt$ can be slid and
isotoped without moving $K$ until $\ttt$ is level -- either a level
arc or a level ``eyeglass'' (the wedge of an arc and a circle). 
Moreover, after $\ttt$ is perturbed slightly (to put $K \cup \ttt$ in
normal position) the graph $K \cup \ttt$ cannot be made thinner by
sliding $\ttt$.
\end{theorem}

For a proof see \cite{GST}.  In fact a similar theorem is true for 
arbitrary genus $2$ spines of $S^{3}$ even when we do not allow edges 
to slide, see \cite{ST2}.  

\bigskip

Theorem \ref{theorem:thin=level} raises the natural question whether a 
similar theorem is true for more than a single tunnel.  In general, 
for $K$ a knot in $S^{3}$, a collection $\ttt_{1}, \ldots, \ttt_{n}$ of 
disjoint properly embedded arcs in $S^{3} - K$ is a system of 
unknotting tunnels if the graph $\Ggg = K \cup (\ttt_{1} \cup \ldots \cup 
\ttt_{n})$ is a Heegaard spine (i.e. the complement of $\eta(\Ggg)$ 
is a handlebody).  

\begin{quest} Suppose $\ttt_{1}, \ldots, \ttt_{n}$ is a system of
unknotting tunnels for a knot $K \subset S^{3}$, in thin position with
respect to the standard height function $h$.  Suppose $\ttt_{1},
\ldots, \ttt_{n}$ are slid and isotoped to minimize the width of $\Ggg =
K \cup (\ttt_{1} \cup \ldots \cup \ttt_{n})$. Is each of the tunnels 
a perturbed level arc?
\end{quest}

Of course many versions of this question are possible, e.  g.
extending it to links or to arbitrary graphs in $S^{3}$.  Even the
case of a pair of tunnels seems difficult; it seems the first order of
business would need to be a generalization of Morimoto's theorem
\cite{Mo} (so essential for the proof of Theorem \cite{Th1}) to
handlebodies of higher genus.

\section{Graphs in other $3$-manifolds}

All our discussion so far revolves around objects (knots, links, 
graphs) in the $3$-sphere, on which we have the standard height function.  
In fact, one can imagine using thin position in many other contexts.  
For example, if $K$ is a knot in an arbitrary closed $3$-manifold $M$, and 
$H_{1} \cup_{P} H_{2}$ is a Heegaard splitting for $M$, one can 
describe the Heegaard splitting as a product structure on the 
complement of spines $\Ss_{1} , \Ss_{2}$ for the respective 
handlebodies $H_{1}, H_{2}$.  That is, $M - (\Ss_{1} \cup 
\Ss_{2}) \cong (P \times (-1, 1))$.  Just as in the applications 
above, one can define the width of $K$ with respect to this 
structure, and try to minimize the width.  In effect, we are 
retrospectively viewing the whole discussion above as the special case 
in which $M$ is $S^{3}$ and the splitting is of genus $0$.  Of course, 
many of the arguments above rely heavily on the fact that $P$ is a 
sphere, so generalizing in this direction has not been particularly 
fruitful.

But there is a remarkable application of thin position that occurs as
a crucial step in Thompson's recognition algorithm for the $3$-sphere
\cite{Th2}.  In this section we will briefly outline how
it arises, and note some related applications to other decision
problems in $3$-manifold topology.

\begin{figure}[tbh]
\centering
\includegraphics[width=.6\textwidth]{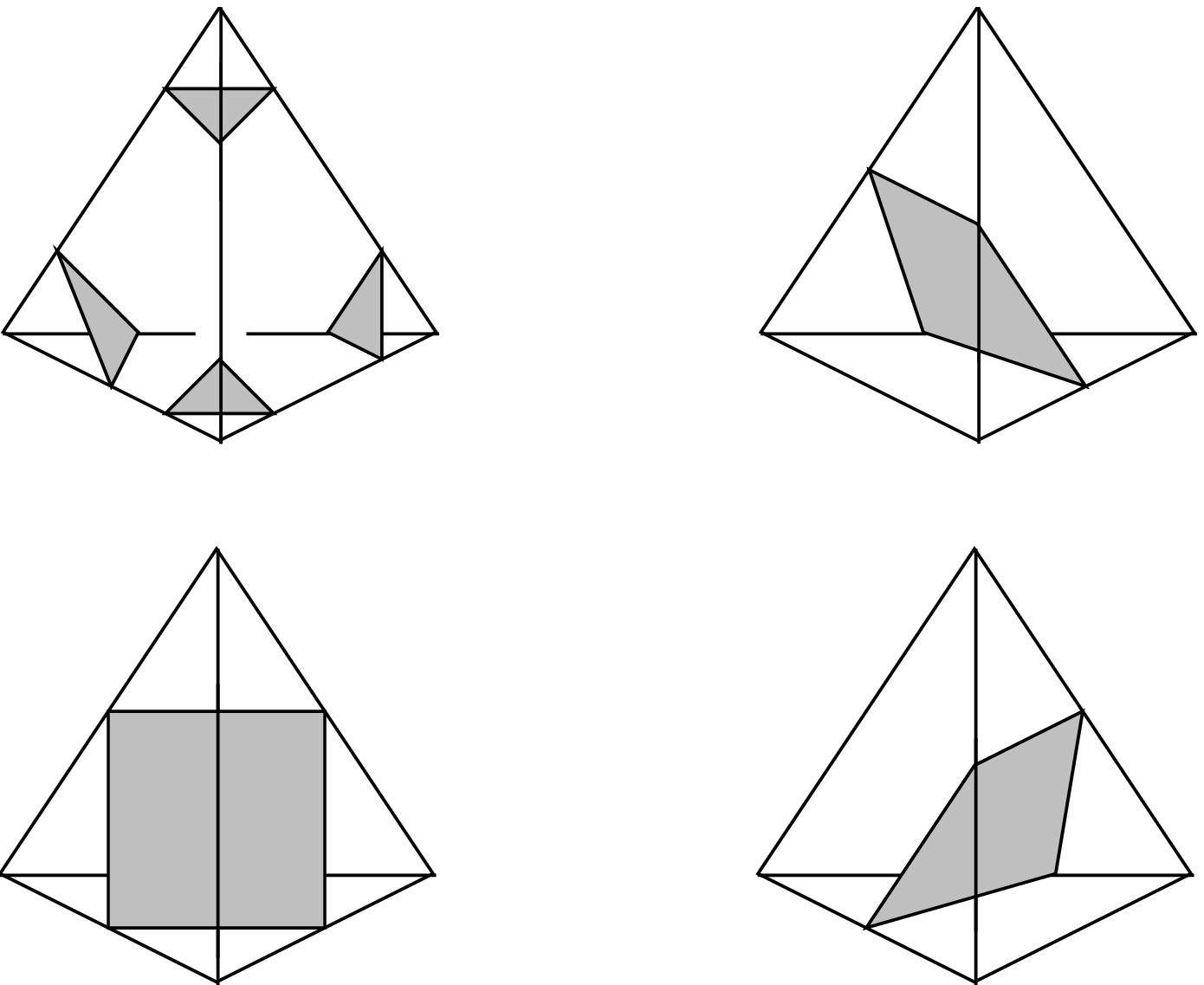}
\caption{} \label{fig:normal}
\end{figure}

Suppose $M$ is a closed $3$-manifold with a given triangulation
$\mathcal{T}$.  Let $\Ggg$ be the $1$-skeleton of the triangulation. 
Recall that a compact surface $F \subset M$ is {\em normal} with
respect to the triangulation if 
\begin{itemize} 
    
    \item $F$ is in general position with
respect to $\mathcal{T}$ (so in particular $F$ intersects $\Ggg$ is a
finite number of points, each on an edge of $\Ggg$)

\item For each $2$-simplex $\Ddd_{2}$ in $\mathcal{T}$, each component 
of $\Ddd_{2} \cap F$ is an arc with its ends on different faces of 
$\Ddd_{2}$.

\item For each $3$-simplex $\Ddd_{3}$ in $\mathcal{T}$, each component 
of $\Ddd_{3} \cap F$ is either a triangle (i. e. parallel to a face 
of $\Ddd_{3}$) or a square (i. e. it is incident to each face in a 
single arc).  See Figure \ref{fig:normal}.

\end{itemize}

It is easy to show that any incompressible surface in $M$ can be
isotoped so that it is normal.  The converse is not true, most
obviously because the link in $M$ of any vertex of $\mathcal{T}$ is a
normal sphere.  On the other hand, if $F$ is a normal surface, then
the complement $F - \Ggg$ is incompressible in $M -
\Ggg$.  That implication is essentially reversible: If $F$ is a
surface so that $F - \Ggg$ is incompressible in $M -
\Ggg$ then $F$ may be isotoped rel $F \cap \Ggg$ so that either
$F$ is normal or it's a sphere that bounds a ball intersecting $\Ggg$ 
in a single unknotted arc.

There is an algorithm to find a maximal collection $\Ss$ of disjoint
non-parallel normal surfaces in $M$; the roots of this algorithm go
back to early work of Kneser, establishing that there are at most a
finite number of connected summands in $M$ \cite{Kn}.  At the very
least, $\Ss$ contains a linking sphere of each vertex, but typically
there are many more.  For example, if an edge of $\mathcal{T}$ is
incident to two distinct vertices, tube together their linking
vertices by a tube along the edge.  Unless $M$ is a specific
$2$-vertex triangulation of $S^{3}$ (cf \cite{JR}) such a sphere is
normal.  In the end, it is possible to show that each component of $M
- \Ss$ is one of three types:

\begin{itemize}
    
    \item a ball containing a single vertex, and bounded by a
    vertex-linking sphere
    
    \item a punctured $3$-ball with more than one boundary component
    
    \item a single further component $M_{0}$, for which $|\bdd 
    M_{0}| = 1$.
    
    \end{itemize}
    
Then $M$ is the $3$-sphere if and only if $M_{0}$ is a $3$-ball, and
$M_{0}$ is algorithmically recognizable among the components of $M -
\Ss$ by the fact that it is the only component that has a single
boundary component and contains no vertex.  So in order to determine
if $M$ is a $3$-sphere, it suffices to find an algorithm to decide if
$M_{0}$ is a $3$-ball.
    
Inside of $M_{0}$ is a proper collection of arcs, $K = \Ggg \cap
M_{0}$ and, because $\Ss - \Ggg$ is incompressible in the complement
of $\Ggg$, we have that the planar surface $\bdd M_{0} - K$ is
incompressible in $M_{0} - K$.  Suppose $M_{0}$ is a $3$-ball and
imagine putting $K$ in thin position with respect to the radial height
function on the ball.  We know immediately that $K$ is also in bridge
position; that is, all the maxima lie above all the minima.  For if
not, consider the thin spheres in $M_{0}$.  We have noted above
(essentially Theorem \ref{theorem:wu}) that some thin sphere $P$ has
the property that $P - K$ is incompressible in $M_{0} - K$; hence $P$
would be a normal sphere in $M_{0}$ not parallel in $M_{0} - K$ to
$\bdd M_{0}$.  But this would contradict the completeness of $\Ss$.

This connection between thin spheres in $M_{0} - K$ and normal spheres
in $M$ prompts this question: what would a thick sphere tell us? 
(Note that there has to be a thick sphere, since if there are only
minima in $K$ then $\bdd M_{0} - K$ would be compressible.)  If $P$
is the (unique) thick sphere, then a maximum of $K$ can be pushed down
to $P$ and a minimum pushed up, but not simultaneously.  Translating
back into how $P$ would appear in the triangulation, it turns out that
it looks just like a normal sphere, except in a single $3$-simplex
$\Ddd_{3}$, where a single component is not a triangle or a square,
but an octagon (cf. Figure \ref{fig:almost}).  Such a surface is called
an {\em almost normal surface}.  Observe that arcs in the $1$-skeleton
of the $3$-simplex can be pushed to the octagon from either side, but
their images there necessarily intersect.  Roughly the same algorithm
that detects normal spheres can be used to detect almost normal 
spheres.  The upshot is this:

\begin{figure}[tbh]
\centering
\includegraphics[width=.3\textwidth]{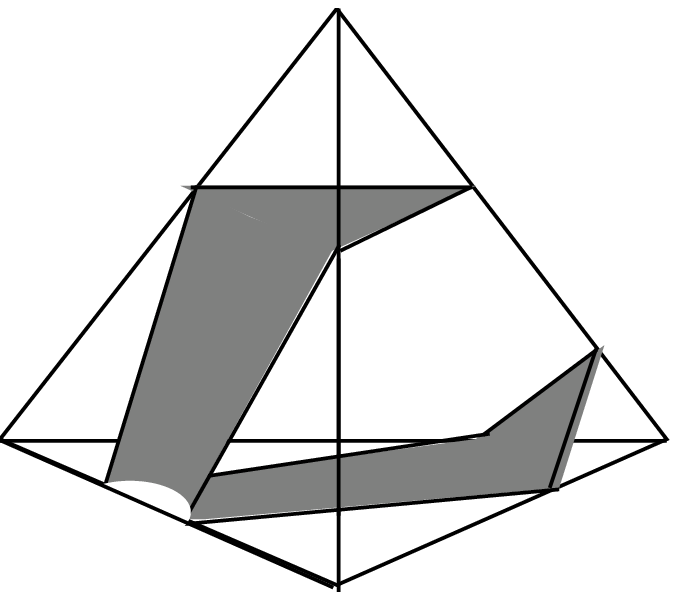}
\caption{} \label{fig:almost}
\end{figure}

\bigskip

{\bf Fact 1:}  If $M$ is the $3$-sphere then there is an almost normal 
sphere in the component $M_{0}$ of $M - \Ss$.

\bigskip

To complete the argument that this is an algorithm, one needs to know
that if there is an almost normal sphere $P$ in the component $M_{0}$
of $M - \Ss$ then $M_{0}$ is not the $3$-ball and so $M$ is not
$S^{3}$.  Observe first of all that since a sub-arc of $\Ggg$ is
parallel to an arc in $P$, it follows that $P - K$ is compressible in
$M_{0} - K$ on the side containing the arc: basically one constructs a
compressing disk by doubling the disk defining the parallism.  This
argument applies on both sides of $P$, so $P$ is compressible in
$M_{0} - K$ in both directions (but compressing disks on opposite
sides necessarily intersect).  Thicken $P$ to a collar $P \times I$,
then maximally compress $(P \times \bdd I) - K$ in $M_{0} - K$; the
result must be an incompressible planar surface (possibly with many
components).  Each component then just comes from a trivial sphere
cutting off a ball intersecting $K$ in an unknotted arc, or it becomes
a normal $2$-sphere, hence a sphere parallel in $M_{0} - K$ to $\bdd
M_{0} - K$.  Filling in all the $3$-balls, then, creates exactly a
copy of $M_{0}$.  That is, $M_{0}$ can be obtained from $P \times I$
by attaching only $2$ and $3$-handles to $P \times \bdd I$.  In
particular, if $M_{0}$ is a homology ball, it's a real ball.  We
conclude:

\bigskip

{\bf Fact 2:} If $M$ is a homology sphere and there is an almost 
normal sphere in the component $M_{0}$ of $M - \Ss$, then $M$ is 
the $3$-sphere.

\bigskip

Since the homology of $M$ is easily calculable, the combination of 
Facts $1$ and $2$ gives the Thompson algorithm for recognizing the 
$3$-sphere.  

\bigskip

Without wandering too far afield, note that the algorithm above is
particularly straightforward if there are few normal spheres.  A
triangulation (broadly defined) of a closed $3$-manifold is {\em
$0$-efficient} if the only normal spheres are vertex linking.  Clearly
such a manifold must be irreducible.  It is a theorem of Jaco and
Rubinstein \cite{JR} that, with just a few specific exceptions, any
triangulation of a closed, orientable, irreducible $3$-manifold can be
modified to be $0$-efficient and such a triangulation has only a
single vertex.  More generally, if $M$ is reducible, there is an
algorithm to decompose $M$ into a connected sum of $3$-manifolds, each
of which either has a $0$-efficient triangulation or is visibly
homeomorphic to $S^{3}, S^{2} \times S^{1}, RP^{3}$ or $L(3,1)$.  Thus
to get an algorithm that precisely describes the connected sum
decomposition of $M$, one need only apply the Thompson algorithm above
in the case in which the triangulation is $0$-efficient.

\bigskip

There are other clever applications of thin position in settings that
go well beyond the scope of this article.  A favorite is \cite{Lac},
where Lackenby shows that the natural combinatorial ideal triangulation of a
punctured torus bundle (with pseudo-Anosov monodromy) coincides with
the natural hyperbolic ideal triangulation.

\end{document}